\documentclass[11pt]{amsart}

\usepackage{pdfpages}
\usepackage{amssymb}
\usepackage{amsmath}
\usepackage{amsthm}
\usepackage{thm-restate}
\usepackage{bbm}
\usepackage{graphicx}
\usepackage{float}
\usepackage{tikz}
\usetikzlibrary{arrows.meta,fit,backgrounds,shapes.geometric}
\graphicspath{ {./pics/} }

\newtheorem{theorem}{Theorem}
\newtheorem{corollary}[theorem]{Corollary}
\newtheorem{lemma}[theorem]{Lemma}
\newtheorem{definition}[theorem]{Definition}
\newtheorem{proposition}[theorem]{Proposition}
\newtheorem*{conjecture}{Conjecture}
\newtheorem*{theorem*}{Theorem}

\numberwithin{theorem}{section} 
\numberwithin{corollary}{section}
\numberwithin{lemma}{section}
\numberwithin{proposition}{section}
\numberwithin{definition}{section}

\begin{document}
\title{The Local H\"older Exponent of Core Entropy}

\author{Jacob Kewarth}
\address{University of Toronto, 40 St. George Street, M5S 2E4, Toronto, Canada}
\email{jacob.kewarth@mail.utoronto.ca}

\begin{abstract}
Here we study the local H\"older exponent of core entropy $h(M)$. For degree $d$ unicritical primitive majors $M$ we show the H\"older exponent equals $h(M) / \log d$ whenever $h(M) > 0$. For multicritical majors we give an explicit counterexample to the relationship.
\end{abstract}

\maketitle


\section{Introduction}

The core entropy $h(f)$ of a postcritically finite (PCF) polynomial $f$ is the topological entropy of $f$ restricted to its Hubbard tree. This quantity has been widely studied, see for instance \cite{DS}, \cite{Ga}, \cite{GT}, \cite{HS}, \cite{Ke}, \cite{Ju}, \cite{Li}, \cite{Ti2}. A common approach is to associate to each PCF polynomial a primitive major, that is, a finite lamination whose leaves represent the critical points of the polynomial. Thus we often instead define core entropy as a map $h(M)$ over the space $PM(d)$ of degree $d$ primitive majors.

A common question about core entropy has been regarding the relationship between the entropy itself and the local H\"older exponent of the entropy. Indeed, \cite{GT} write ``we believe the local H\"older exponent of the entropy function is related to the value of the entropy'', and \cite{HS} explicitly state as a conjecture ``the H\"older exponent is equal (up to a scaling constant) to core entropy itself''. Such a phenomenon has been observed for other dynamical settings, see, for example: \cite{CT}, \cite{DM}, and \cite{Ti3}.

More precisely, if we equip $PM(d)$ with the Thurston metric, we ask if the local H\"older exponent (see section \ref{pmd} for definitions) of core entropy $h$ satisfies:
\begin{conjecture}
Fix $d \geqslant 2$. For all $M \in PM(d)$, if $h(M) > 0$ then:
\[
\alpha_{loc}(h, M) = \frac{h(M)}{\log d}.
\]
\end{conjecture}

Note that we restrict to $h(M) > 0$ as when $h(M) = 0$ it has been observed that H\"older continuity fails even though the function is globally continuous (see \cite{Ti2}).

\cite{Ju} estimated the asymptotic behaviour of core entropy along certain sequences of rational angles, supporting the conjecture. For a large class of quadratics the lower bound on the exponent is shown in \cite{HS} (along with symbolic estimates on the upper bound in terms of kneading sequences). However even for quadratics the conjecture has remained open. The goal of this paper is to settle the conjecture.  In particular, we wish to prove it for unicriticals and disprove it in general.

As in the approach of \cite{HS} we begin by defining a symbolic core entropy function $h(w)$ on spaces of sequences $\{0, 1, \dots, d-1 \}^{\mathbb{N}}$. Here we rely less on Hubbard trees and more on the labeled wedge techniques of \cite{Ti2} and \cite{GT} which we recall in section \ref{background}. In this theory we define an infinite graph for each sequence $w$ (or primitive major $M$) and then relate the entropy to the growth rate of this graph which occurs as the inverse of the smallest root of the spectral determinant of the graph. In section \ref{firstreturns} we develop an alternative definition. We show that one can extract the growth rate of the graph by considering only the growth rate of paths that visit a finite portion of the graph, and that the size of this region is locally uniform in the entropy. This leads to the introduction of our first return matrix $E_M(z)$ (see definition \ref{firstreturnmatrix}), a finite matrix whose dimension is locally uniform with respect to the entropy and whose leading eigenvalue determines the entropy of $w$. 

Section \ref{pathcounting} uses this construction to prove proposition \ref{4prop} which says that the entropy root of the spectral determinant is simple. Section \ref{entropygaps} develops estimates needed for the upper bound on the H\"older exponent.

With this technology established section \ref{secsymbol} proves that for the metric $D(w, w') = d^{-N(w, w')}$, where $N(w, w')$ is the first index where the sequences differ:

\begin{theorem}
\label{symbolicthm}
Fix $d \geqslant 2$. For all $w \in \{0, 1, \dots, d-1\}^{\mathbb{N}}$, if $h(w) > 0$ then:
\[
\alpha_{loc}(h, w) = \frac{h(w)}{\log d}.
\]
\end{theorem}

The work up till this point is purely symbolic dynamics. Section \ref{anglelemmas} gives us tools to convert these results into information about primitive majors.

Then in section \ref{mainproofsection} we show the main result of the paper:

\begin{theorem}
\label{thm}
Fix $d \geqslant 2$. For all $M \in \Pi(d)$, if $h(M) > 0$ then:
\[
\alpha_{loc}(h \mid_{\Pi(d)}, M) = \frac{h(M)}{\log d}.
\]
\end{theorem}

Finally the paper ends with section \ref{counterexmaplesection} which gives an explicit counterexample to the general conjecture.

\begin{theorem}
There exists $M \in PM(3)$ such that:
\[
\alpha_{loc}(h, M) < \frac{h(M)}{\log 3}.
\]
\end{theorem}

\medskip
\noindent\textbf{AI Disclosure.}
The author used ChatGPT (OpenAI) to assist with mathematical
discussions, generating figures, and editing during the preparation of this paper.
The final draft was entirely written by the author who takes full responsibility for its contents.

\section{Background}
\label{background}

\subsection{Graphs with Bounded Cycles}

Let $\Gamma$ be an infinite graph with bounded outgoing edges. Assume also that $\Gamma$ has bounded cycles, that is, that the number of simple cycles of length $\leqslant n$ is finite for all $n$.

Define the growth rate of $\Gamma$ to be:

\[
\rho(\Gamma) = \limsup_{n \rightarrow \infty} C_\Gamma(n)^{1/n}
\]

where $C_\Gamma(n)$ is the number of marked closed paths of length $n$. By marked we mean that we are counting paths that visit the same vertices as different if their starting vertex is different. Define the inverse growth:

\[
r(\Gamma) = \frac{1}{\rho(\Gamma)}.
\]

The spectral determinant is:

\[
P_\Gamma(z) = \sum_{\gamma} (-1)^{c(\gamma)}z^{l(\gamma)}
\]

where the sum is taken over all simple multi-cycles in $\Gamma$. $l(\gamma)$ denotes the length of the multi-cycle and $c(\gamma)$ denotes the number of components.

Also let:

\[
\sigma(\Gamma) = \limsup_{n \rightarrow \infty} \text{(\# of multi-cycles of length $n$)}^{1/n}
\]

The following is \cite[theorem 4.3]{Ti2}:
\begin{lemma}
If $\sigma(\Gamma) \leqslant 1$ then $P_\Gamma(z)$ is holomorphic in the unit disc. Moreover if $r(\Gamma) < 1$ then $P_\Gamma(z)$ is non-zero on $|z| < r(\Gamma)$ and $P_\Gamma(r(\Gamma)) = 0$.
\end{lemma}

Since the graph $\Gamma$ is infinite its corresponding adjacency matrix $A$ is infinite. Still, because the number of outgoing edges is bounded and the graph has bounded-cycles, many facts of finite linear algebra are still true. In particular, the following is a well defined finite number for each $n$:

\[
(A^n)_{u, v} := \text{The number of paths of length $n$ from $u$ to $v$}.
\]

Hence we may define:

\[
Tr(A^n) = \sum_{v \in V} (A^n)_{v, v}.
\]

For each $n$ only a finite number of the above terms in the sum are non-zero so we avoid any convergence issues.

Note by definition we have $C_\Gamma(n) = Tr(A^n)$.

Then by \cite[lemma 4.4]{Ti2} we have:

\[
\frac{1}{P_\Gamma(z)} = \exp(\sum_n \frac{Tr(A^n)}{n}z^n)
\]

For each $n \in \mathbb{N}$ and each vertex $u$ of $\Gamma$ we define:
\[
B_\Gamma(u, n) = \# \{ \text{Paths of length $n$ in $\Gamma$ starting at $u$} \}.
\]

\subsection{Wedges}

Our goal here is to define a broad family of graphs satisfying the conditions from before.

Fix an integer $d$ and an integer $1 \leqslant s \leqslant d-1$. We call a wedge of size $s$ a set:

\[
W_s = \Big\{ \{ (k, i) , (l, j) \} \mid \begin{gathered} k, l \in \{1, 2, \dots, s \}, \\
1 \leqslant i \leqslant j, (k, i) \neq (l, j)
\end{gathered}
 \Big\}.
\]

Then our $W_s$ will be the vertex set for each graph studied while the edge set will vary. For a vertex $\{ (k, i), (l, j)\}$ we call $i$ the height of the vertex and $j$ the width.

A labeling $L$ of $W_s$ is a map taking a vertex $\{ (k, i), (l, j)\}$ and assigning it either $\emptyset$ or $(\alpha_1, \alpha_2, \dots, \alpha_r)$ with $\alpha_i$ pairwise distinct elements in $\{1, 2, \dots, s \}$. We call the vertex non-separated in the first case and separated in the second case and refer to $(\alpha_1, \dots, \alpha_r)$ as the separation vector of $\{ (k, i), (l, j) \}$.

Given a labeling $L$ of $W_s$ we can define a graph $\Gamma$ as follows. Each element of $W_s$ is a vertex. For each vertex $v = \{ (k, i), (l, j) \}$ we consider its label. If the label is $\emptyset$ then there is one edge:

\[
\{ (k, i), (l, j) \} \rightarrow \{ (k, i+1), (l, j+1) \}.
\]
If the separation vector is $(\alpha_1, \dots, \alpha_r)$ then there are $r+1$ edges:

\begin{align*}
\{ (k, i), (l, j) \} &\rightarrow \{ (k, i+1), (\alpha_1, 1) \} \\
&\rightarrow \{ (\alpha_1, 1), (\alpha_2, 1)\} \\
&\rightarrow \{ (\alpha_2, 1), (\alpha_3, 1) \} \\
& \dots\\
&\rightarrow \{(\alpha_r, 1), (l, j+1) \}.
\end{align*}

Then clearly any such graph $\Gamma$ has finitely many outgoing edges. Moreover, according to \cite[corollary 4.2]{GT} we have that the graph has bounded cycles.

If a vertex is separated then we denote the edge $\{ (k, i), (l, j) \} \rightarrow \{ (\alpha_r, 1), (l, j+1)\}$ the forward edge and call all other edges backwards edges.

In the special case when $s = 1$, an element $\{ (k, i), (l, j) \}$ of our wedge $W_s$ has $k = l = 1$. In this case we omit both elements and view $W_1$ as a set of tuples $(i, j)$ with $i < j$ instead. For convenience we will often treat $(i, j)$ as being unordered to avoid having to track which index is larger.

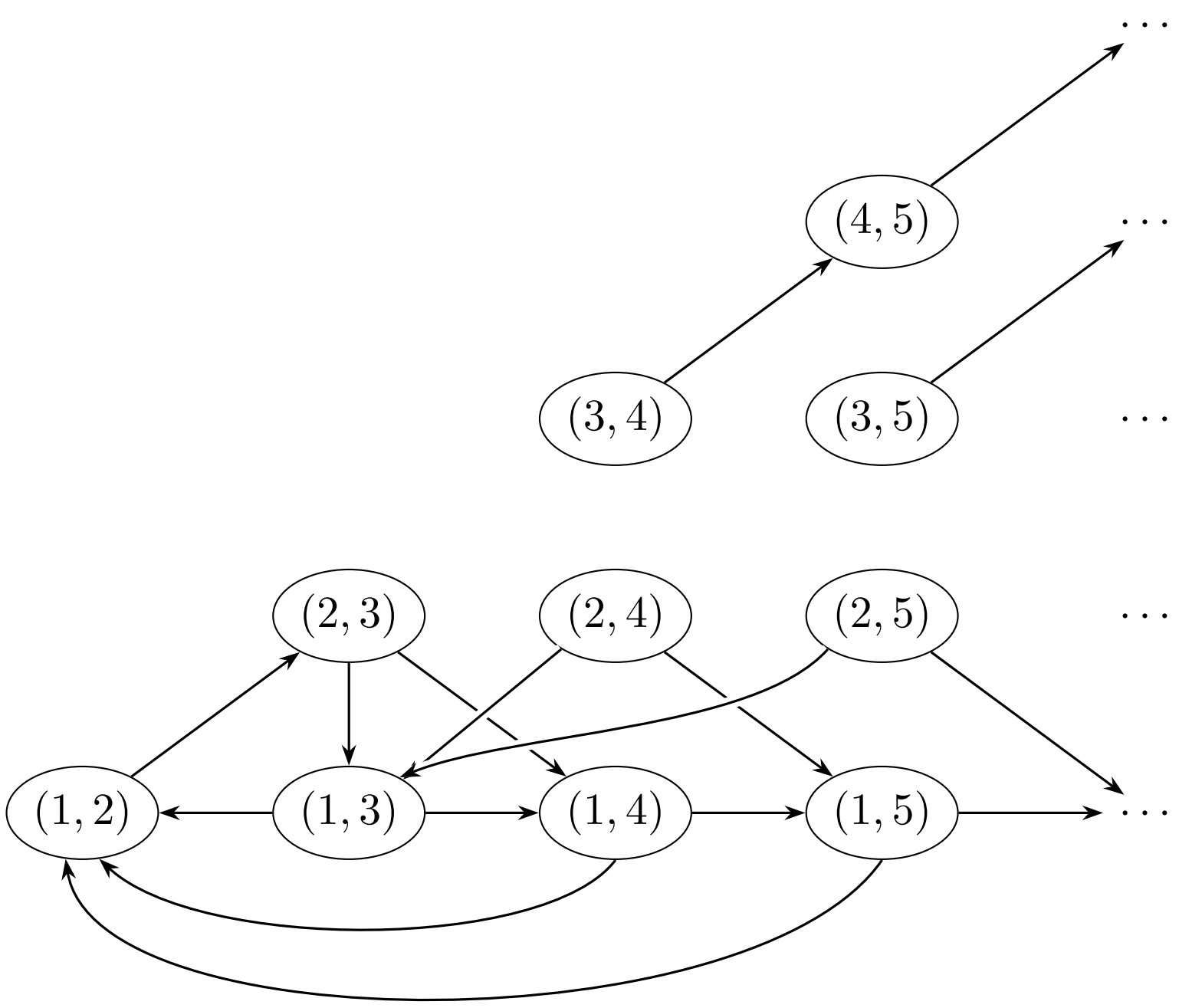
\begin{figure}[htbp]
\centering
\begin{tikzpicture}[
  >={Stealth[length=1.8mm,width=1.25mm]},
  vertex/.style={
    ellipse,draw,minimum width=13mm,
    minimum height=8mm,inner sep=1pt
  },
  every edge/.style={draw,->,semithick},
  font=\normalsize
]
  \node[vertex] (12) at (0,0) {$(1,2)$};
  \node[vertex] (13) at (2.3,0) {$(1,3)$};
  \node[vertex] (14) at (4.6,0) {$(1,4)$};
  \node[vertex] (15) at (6.9,0) {$(1,5)$};

  \node[vertex] (23) at (2.3,1.7) {$(2,3)$};
  \node[vertex] (24) at (4.6,1.7) {$(2,4)$};
  \node[vertex] (25) at (6.9,1.7) {$(2,5)$};

  \node[vertex] (34) at (4.6,3.4) {$(3,4)$};
  \node[vertex] (35) at (6.9,3.4) {$(3,5)$};

  \node[vertex] (45) at (6.9,5.1) {$(4,5)$};

  \node (16) at (9.2,0) {$\dots$};
  \node at (9.2,1.7) {$\dots$};
  \node at (9.2,3.4) {$\dots$};
  \node (46) at (9.2,5.1) {$\dots$};
  \node (56) at (9.2,6.8) {$\dots$};

  \path
    (12) edge (23)
    (13) edge (12)
    (13) edge (14)
    (14) edge (15)
    (23) edge (13)
    (23) edge (14)
    (24) edge (15)
    (34) edge (45)
    (15) edge (16)
    (25) edge (16)
    (35) edge (46)
    (45) edge (56);

  \draw[->,semithick] (14.south)
    .. controls (4,-1.2) and (1,-1.2)
    .. (12.290);

  \draw[->,semithick] (15.south)
    .. controls (5.8,-2) and (0.1,-2)
    .. (12.250);

  \draw[->,semithick,
    preaction={draw=white,line width=3pt}]
    (24.south west) -- (13.north east);

  \draw[->,semithick,
    preaction={draw=white,line width=3pt}]
    (25.south west)
    .. controls (5.8,0.7) and (3.7,0.7)
    .. (13.35);
\end{tikzpicture}
\caption{The labeled wedge corresponding to the kneading sequence $w=001111\dots$.}
\label{fig:labelled-wedge-example}
\end{figure}

\subsection{Kneading Sequences}

Fix $d \geqslant 2$. Then an unstarred kneading sequence is a sequence $w \in \{0, 1, \dots, d-1 \}^{\mathbb{N}}$. A starred kneading sequence (or $\star$ kneading sequence) is a sequence of the form $w = (w_1w_2\dots w_{p-1}\star)^\infty$ for $w_1, \dots, w_{p-1} \in \{0, 1, \dots, d-1 \}$. Note here a $\star$ kneading sequence is necessarily periodic with period equal to the first index which is a $\star$. Given a kneading sequence $w = w_1w_2\dots$, whether starred or unstarred we can define a labeling as follows. Given a pair $i \neq j$, if $w_i$ and $w_j$ are not $\star$ then we label $(i, j)$ separated if $w_i \neq w_j$ and unseparated if $w_i = w_j$. If one (or both) of $w_i$ and $w_j$ is a $\star$, we label the vertex unseparated. In this way, each kneading sequence $w$ induces a labelled wedge graph $\Gamma_w$ (see figure \ref{fig:labelled-wedge-example}). Then we can define the core entropy of our kneading sequence by:
\[
h(w) = \log \max(1, \rho(\Gamma_w)).
\]
Note that if $h(w) > 0$ then we have $h(w) = -\log r(\Gamma_w)$.

\cite[lemma 7.3]{Ti2} gives a useful way to compare entropies of starred and unstarred kneading sequences:
\begin{lemma}
\label{starcomparison}
Given a starred kneading sequence $w = (w_1\dots w_{p-1} \star)^\infty$, let $\tau \in \{0, 1, \dots, d-1 \}$ be arbitrary. Then for $u = (w_1\dots w_{p-1}\tau)^\infty$ we have:
\[
h(w) = h(u).
\]
\end{lemma}
Because of this we often avoid discussion of starred kneading sequences during symbolic proofs.

Given two kneading sequences $u$ and $w$, we can define $N(u, w)$ to be the first index where their symbols differ. Then we can define their distance as:
\[
D(u, w) = d^{-N(u, w)}.
\]
Then for unstarred kneading sequences, \cite[lemma 4.4]{GT} gives continuity:
\begin{lemma}
\label{symboliccontinuity}
If $w_n, w$ are unstarred kneading sequences and $w_n \rightarrow w$ then $h(w_n) \rightarrow h(w)$.
\end{lemma}

\subsection{Finite Graph Model}

Here we describe an alternative finite graph that we can associate to a kneading sequence from which we can study its core entropy. As the finite graph implies a finite adjacency matrix this is sometimes preferred to the labeled wedge construction from before, however, a major drawback to this approach is that the number of vertices of the graph is not locally bounded.

Consider $w$ a periodic or preperiodic unstarred kneading sequence, say with preperiod $q \geqslant 0$ and period $p \geqslant 1$. Note for this construction we do not assume the preperiod or period is minimal, thus strictly speaking this construction depends both on our kneading sequence and on our choice of $q$ and $p$. In the event we do not specify the preperiod and period in our construction they should be taken to be minimal. If we have $w = (w_1...w_p)^\infty$ then we construct the graph with preperiod $0$ and period $p$.

Define a vertex set:
\[
V = \big\{ (i, j) \mid 1 \leqslant i < j \leqslant q+p \big\}.
\]
Note this set is nonempty when $p+q \geqslant 2$. If $p+q = 1$ then our sequence must be a constant sequence and hence have zero entropy. Thus we shall assume $p+q \geqslant 2$.

We label each vertex as before. If a vertex $(i, j)$ is non-separated then it has one edge:
\[
(i, j) \rightarrow \begin{cases}
(i+1, j+1) & j < q+p \\
(i+1, q+1) & j = q+p
\end{cases}
\]
If the vertex is separated, we have two edges:
\begin{align*}
(i, j) &\rightarrow (1, i+1) \\
(i, j) &\rightarrow \begin{cases}
(1, j+1) & j < q+p \\
(1, q+1) & j = q+p
\end{cases}
\end{align*}
Note that in the special case that the target in any of these cases would be diagonal, then we say there is no edge. 

Then we call this graph $G_w$. We can extend this to starred kneading sequences $w$ by replacing them with an unstarred kneading sequence as described in the previous section. Then we can repeat this construction with $q = 0$ and $p$ the same period as for the $\star$ kneading sequence.

According to \cite[proposition 7.2]{Ti2} we have:
\begin{proposition}
\label{tiozzo5.12}
Let $w$ be a periodic or preperiodic kneading sequence. Let $\lambda$ be the leading eigenvalue of the adjacency matrix of the graph $G_w$. Then:
\[
h(w) = \log \max(1, \lambda).
\]
In particular, if $h(w) > 0$ then $\lambda = r(\Gamma_w)^{-1}$.
\end{proposition}

\subsection{PM(d)}\label{pmd}
Fix $d \geqslant 2$. Then a collection $M = \{L_1, \dots, L_k \}$ of ideal lines or polygons (called leaves) on $\overline{\mathbb{D}}$ is said to be a degree $d$ primitive major if and only if:
\begin{enumerate}
\item The leaves $L_j$ are pairwise disjoint,
\item $\sum_{j=1}^k \big( |L_j \cap \partial \mathbb{D}| - 1 \big)= d-1$,
\item for each $j$, the endpoints $L_j \cap \partial \mathbb{D}$ are identified by $z \mapsto z^d$.
\end{enumerate}

Then we call $PM(d)$ the collection of all such primitive majors.

Each primitive major $M$ induces a quotient graph $\gamma(M)$ by taking $\partial \mathbb{D}$ and identifying each leaf of $M$. Then the path metric on the unit circle gives a metric on $\gamma(M)$ which we can pull back to get a pseudometric on $\partial \mathbb{D}$ which we will denote by $d_M$. Then we can define a metric $md$ on $PM(d)$ by:
\[
md(M, M') = \sup_{x, y \in \partial \mathbb{D}} |d_M(x, y) - d_{M'}(x, y)|.
\]
Here we identify $\partial \mathbb{D}$ with $\mathbb{R} / \mathbb{Z}$. In particular, when discussing the path metric on the unit circle we assume it is normalized so the entire circle has length $1$. See \cite{Th+} or \cite{GT} for more details.

Given a primitive major $M$, we may associate to it a labelled wedge $\Gamma_M$ by looking at iterates of the leaves of $M$ (see \cite[section 4.4]{GT}) and in this way define core entropy as a map on $PM(d)$:
\[
h(M) = \log \max(1, \rho(\Gamma_M)).
\]
Under the metric $md$, we have that $h(M)$ is continuous \cite{GT}. Using this metric we can further define the local H\"older exponent by:
\[
\alpha_{loc}(h, M) = \sup\left\{ \gamma>0 \mid
\begin{array}{l}
\exists C,\epsilon>0, \text{ if } md(M, M') <\epsilon \text{ then } \\
|h(M) - h(M') |\leqslant C md(M, M')^\gamma
\end{array}
\right\}.
\]

\begin{figure}[htbp]
\centering
\begin{tikzpicture}[font=\normalsize]
  \def\R{1.45}

  \begin{scope}
    \draw[gray] (0,0) circle[radius=\R cm];
    \foreach \a in {45,135,225,315} {
      \draw[thick] (\a:\R)
        arc[start angle={\a-90},
            end angle={\a-180},radius=\R cm];
    }
    \node at (0,-2) {$\Pi(4)$};
  \end{scope}

  \begin{scope}[xshift=4.2cm]
    \draw[gray] (0,0) circle[radius=\R cm];
    \foreach \a in {30,120} {
      \draw[thick] (\a:\R)
        arc[start angle={\a-90},
            end angle={\a-180},radius=\R cm];
    }
    \draw[thick] (210:\R) -- (30:\R);
    \draw[thick] (240:\R)
      arc[start angle=150,end angle=60,radius=\R cm];
    \node at (0,-2) {$\Pi(3,2)$};
  \end{scope}

  \begin{scope}[xshift=8.4cm]
    \draw[gray] (0,0) circle[radius=\R cm];
    \foreach \a in {10,130,250} {
      \draw[thick] (\a:\R)
        arc[start angle={\a-90},
            end angle={\a-180},radius=\R cm];
    }
    \node at (0,-2) {$\Pi(2,2,2)$};
  \end{scope}
\end{tikzpicture}
\caption{The various strata of $PM(4)$.}
\label{fig:quartic-strata}
\end{figure}
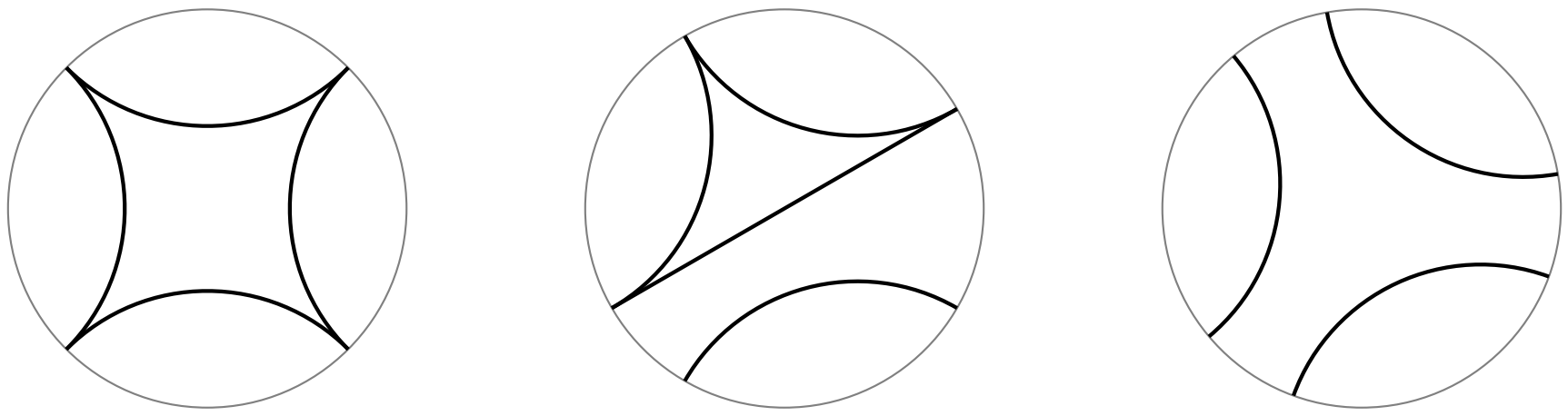

Since the number of leaves in $PM(d)$ can vary it is often common to divide $PM(d)$ into various strata. Specifically, if we have numbers $2 \leqslant s_j \leqslant d$ such that $\sum_{j} (s_j - 1) = d-1$, then we define $\Pi(s_1, \dots, s_k)$ to be the primitive majors of $PM(d)$ of the form $M = \{ L_1, \dots, L_k \}$ such that $|L_j \cap \partial \mathbb{D}| = s_j$ (see figure \ref{fig:quartic-strata}). Of particular importance is the unicritical stratum $\Pi(d)$ which consists of a single leaf with $d$ endpoints. Since a choice of endpoints uniquely determines the leaf, we can represent all such majors by:
\[
M_\theta = \{L_\theta \}, \quad L_\theta = \left\{ \frac{\theta}{d}, \frac{\theta +1}{d}, \dots, \frac{\theta + d - 1}{d} \right\}.
\]
Thus we can parameterize $\Pi(d)$ in terms of $\theta \in \mathbb{R} / \mathbb{Z}$. In fact, when restricted to $\Pi(d)$, $md$ is essentially the same as angular distance:
\begin{lemma}
\label{angularmetric}
For any $d \geqslant 2$ and any $\theta, \theta' \in \mathbb{R} / \mathbb{Z}$ we have:
\[
md(M_\theta, M_{\theta'}) = \frac{2}{d} |\theta - \theta'|.
\]
\end{lemma}

\begin{proof}
Fix $\theta, \theta'$ such that $0 < \theta' - \theta = |\theta' - \theta| \leqslant 1/2$. Let $\epsilon = (\theta' - \theta) / d$. Since rotating the leaf of $M_\theta$ by $\epsilon$ counterclockwise maps $M_\theta$ to $M_{\theta'}$, we have:
\[
d_{M_{\theta'}}(x, y) = d_{M_\theta}(x-\epsilon, y-\epsilon).
\]
Adding and subtracting $d_{M_\theta}(x, y-\epsilon)$ and using the triangle and reverse triangle inequalities gives:
\begin{align*}
|d_{M_\theta}(x-\epsilon, y-\epsilon) - d_{M_\theta}(x, y)| &\leqslant |d_{M_\theta}(x-\epsilon, y-\epsilon) - d_{M_\theta}(x, y-\epsilon)| \\
&+ |d_{M_\theta}(x, y-\epsilon) - d_{M_\theta}(x, y)| \\
&\leqslant d_{M_\theta}(x-\epsilon, x) + d_{M_\theta}(y-\epsilon, y) \\
&\leqslant 2\epsilon.
\end{align*}
Thus we have:
\[
|d_{M_{\theta'}}(x, y) - d_{M_\theta}(x, y)| \leqslant 2\epsilon.
\]
Taking the supremum over all $x, y$ gives:
\[
md(M_\theta, M_{\theta'}) \leqslant 2\epsilon.
\]
Now consider two endpoints $u = \theta/d$ and $v = (\theta + 1) / d$ of $M_\theta$. Since they are identified in the quotient $\gamma(M_\theta)$, we have $d_{M_\theta}(u, v) = 0$. By definition, the distance from $u$ or $v$ to the nearest endpoint of the leaf of $M_{\theta'}$ on $\partial \mathbb{D}$ is $\min(\epsilon, 1/d - \epsilon)$ which is just $\epsilon$ as $\epsilon \leqslant \frac{1}{2d}$. In the quotient $\gamma(M_{\theta'})$ these endpoints are all identified to a single vertex. Since any path in $\gamma(M_{\theta'})$ from $u$ to $v$ must pass through this vertex, we have:
\[
d_{M_{\theta'}}(u, v) = 2\epsilon.
\]
Thus:
\[
md(M_{\theta}, M_{\theta'}) \geqslant 2\epsilon.
\]
Putting these inequalities together and substituting back in the definition of $\epsilon$ we obtain:
\[
md(M_\theta, M_{\theta'}) = \frac{2}{d} |\theta - \theta'|,
\]
as wanted.
\end{proof}

\begin{corollary}
For $d \geqslant 2$, for $M = M_\theta \in \Pi(d)$, if $h(M) > 0$ then:
\[
\alpha_{loc}(h \mid_{\Pi(d)}, M_\theta) = \sup\left\{ \gamma>0 \mid
\begin{array}{l}
\exists C,\epsilon>0, \text{ if } |\theta - \theta'| <\epsilon \text{ then } \\
|h(M_\theta) - h(M_{\theta'}) |\leqslant C |\theta - \theta'|^\gamma
\end{array}
\right\}.
\]
\end{corollary}

Finally, each $M_\theta \in \Pi(d)$ with $\theta < 1$ partitions the unit circle into $d$ pieces which we can label $I_j = (\frac{\theta + j}{d}, \frac{\theta + j + 1}{d})$ for $j = 0, 1, \dots, d-1$. Note here that, unless $\theta = 0$, we have that $0$ belongs to $I_{d-1}$. We can define the kneading sequence of $M_\theta$, or just the kneading sequence of $\theta$, by looking at the orbit of $d^{n-1}\theta \mod 1$ and tracking which interval each iterate lies in, assigning a $\star$ if an iterate hits an endpoint. Denote the kneading sequence by $K(M_\theta)$ or $K(\theta)$. Using this kneading sequence we can define a labeled wedge and so each unicritical primitive major $M$ induces a labeled wedge $\Gamma_M = \Gamma_{K(\theta)}$ of size $1$. 

Since this labeling rule agrees with the one in \cite{GT}, we have that:
\[
h(M_\theta) = h(K(\theta)) = h(\Gamma_{K(\theta)}).
\]

\section{First Returns}
\label{firstreturns}

Here we restrict to the $s = 1$ case. Consider a labeled wedge and let $\Gamma$ denote its induced graph. Then the vertex set consists of tuples $(i, j)$ with $i < j$. For each integer $M \geqslant 2$ we can consider a subset of the vertices:
\[
S_M = \{ (1, j) \mid 2 \leqslant j \leqslant M \}.
\]
For each $M$ removing the vertices $S_M$ from $\Gamma$ yields a new graph which we denote by $\Gamma_M$. The main object of study in this section is:

\begin{definition}\label{firstreturnmatrix}
For each pair $u, v \in S_M$ (allowing $u = v$), let:
\[
e_{uv, M}(n) = \#\{\text{paths $x_0x_1x_2\dots x_n$ in $\Gamma$ with $x_0 = u$, $x_n = v$, and $x_j \notin S_M$ for $j \neq 0, n$ } \}.
\]
Then we define:
\[
E_{uv, M}(z) = \sum_{n=1}^\infty e_{uv}(n)z^n, \quad E_M(z) = (E_{uv, M}(z))_{u, v \in S_M}.
\]
\end{definition}

In the case that our $M$ is fixed we will often omit it and simply reference $e_{uv}$, $E_{uv}(z)$, and $E(z)$ instead.

The paths counted in the definition of $e_{uv}$ are called first return paths and $E_{uv}$ is the first return series. The matrix $E$ is the first return matrix. Notice that since $S_M$ is a finite vertex set, $E$ is a finite matrix.

The goal of this section is to collect facts about our first return matrix. In particular, we prove $\det(I-E(z))$ is related to the spectral determinant before, and hence sees the inverse growth rate of our graph. Moreover we give an estimate on how the size of $E(z)$ changes as we change our graph $\Gamma$ which will be useful for our H\"older estimate later.

\begin{lemma}
\label{countlemma1}
Let $\Gamma$ be the graph of some labeled wedge of size $1$ and fix $M$. Then for each $u \in \Gamma_M$ and each $n \in \mathbb{N}$ we have:
\[
B_{\Gamma_M}(u, n) \leqslant 2 \cdot 2^{\frac{n}{M}}.
\]
Moreover, for every $u \in S_M$ and every $n \in \mathbb{N}$ we have:
\[
\sum_{v \in S_M} e_{uv}(n) \leqslant 8 \cdot 2^{\frac{n}{M}}.
\]

\end{lemma}

\begin{proof}
Consider a separated pair $(i, j)$ not in $S_M$. Then in the graph $\Gamma_M$ this pair has two outgoing edges if $(1, i+1) \notin S_M$, that is, when $i \geqslant M$. Each successor has height $1$, thus one must travel a path of length at least $M-1$ from such a vertex to find another successor with two outgoing edges in $\Gamma_M$. This means that the number of paths of length $n$ starting at an arbitrary vertex $u \in \Gamma_M$ is at most $2\cdot 2^{\frac{n}{M}}$ proving the first claim.

For the second claim consider a first return path of length $n$. It consists of a starting vertex $u$ and an ending vertex $v$, both in $S_M$, along with a path of length $n-2$ in $\Gamma_M$. For each fixed path of length $n-2$ there are at most two choices of the starting and ending vertex. By the above there are at most $2 \cdot 2^{\frac{n-2}{M}}$ choices of paths of length $n-2$ for each starting vertex, hence we get an upper bound of $8 \cdot 2^{\frac{n-2}{M}} < 8 \cdot 2^{\frac{n}{M}}$ as claimed.
\end{proof}

\begin{corollary}
\label{nozerocor}
Let $\Gamma$ be the graph of some labeled wedge of size $1$ and fix $M$. Then each element of $E_M$ is holomorphic on $|z| < 2^{\frac{-1}{M}}$ and $P_{\Gamma_M}(z)$ has no zeros in there.
\end{corollary}

\begin{proof}
The first claim is a routine radius of convergence calculation.

For the second claim note that any vertex on a closed path of length $n$ in $\Gamma_M$ can have width at most $2n+2$. There are at most $2(n+1)^2$ such vertices, hence by lemma \ref{countlemma1}:
\[
C_{\Gamma_M}(n) \leqslant 4(n+1)^2 2^{\frac{n}{M}}.
\]
Thus we have that the series:
\[
\sum_{n=1}^\infty \frac{C_{\Gamma_M}(n)}{n}z^n
\]
has radius of convergence at least $2^{\frac{-1}{M}}$. This means it is holomorphic in this disk, hence $\frac{1}{P_{\Gamma_M}(z)}$ is also holomorphic in the disk and hence the denominator is never $0$.
\end{proof}

\begin{lemma}
\label{unifconv}
Let $\Gamma$ and $\Gamma'$ be the graphs corresponding to two labelled wedges of size $1$. Fix $M$ and assume the labelings of $\Gamma$ and $\Gamma'$ agree at every vertex of width less than $N$, for some $N > M$. Let $E$ and $E'$ denote their respective first return matrices. Then the coefficients of their first return series agree for all $n \leqslant N-M-1$. Moreover, for any $0 < a < 2^{\frac{-1}{M}}$:
\[
\sup_{|z| \leqslant a} ||E(z) - E'(z)||_\infty \leqslant \frac{16}{1- a 2^{\frac{1}{M}}} a^{N-M} 2^{\frac{N-M}{M}}.
\]
where $||E||_\infty = \max_u \sum_{v} |E_{uv}|$ is the usual infinity norm.
\end{lemma}

\begin{proof}
Note that following an edge in $\Gamma$ or $\Gamma'$ can only increase the largest index of a vertex by $1$. Hence, since any vertex in $S_M$ has largest index $\leqslant M$ it follows that any vertex on a path of length $n$ starting at a vertex in $S_M$ has largest index at most $n+M$. Thus if $n \leqslant N-M-1$, the indices of any such vertex are $< N$ meaning their labels in the graphs $\Gamma$ and $\Gamma'$ agree. Since the labels agree the path counts in the first return series will agree proving the first claim.

For the second claim, let $e_{uv}(n)$ and $e'_{uv}(n)$ denote the number of first return paths for the graphs $\Gamma$ and $\Gamma'$ respectively. For each $u \in S_M$, we know by lemma \ref{countlemma1} that:
\[
\sum_{v \in S_M} |e_{uv}(n) - e'_{uv}(n)| \leqslant 16 \cdot 2^{\frac{n}{M}}.
\] 
Also, we know $e_{uv}(n) = e'_{uv}(n)$ for all $n \leqslant N-M-1$. Thus, for any $|z| < a$:
\begin{align*}
\sum_{v \in S_M} |E_{uv}(z) - E'_{uv}(z)| &\leqslant \sum_{n=N-M}^\infty |z|^n \sum_{v \in S_M} |e_{uv}(n) - e'_{uv}(n)| \\
&\leqslant 16 \sum_{n=N-M}^\infty (a2^{\frac{1}{M}})^n \\
&= \frac{16}{1- a 2^{\frac{1}{M}}} a^{N-M} 2^{\frac{N-M}{M}}.
\end{align*}
Taking the maximum over all $u$ and the supremum over $|z| \leqslant a$ then yields:
\[
\sup_{|z| \leqslant a} ||E(z) - E'(z)||_\infty \leqslant \frac{16}{1- a 2^{\frac{1}{M}}} a^{N-M} 2^{\frac{N-M}{M}}
\]
as wanted.
\end{proof}

Finally with these estimates we are now able to relate the determinant of our first return matrix to the spectral determinants $P_\Gamma$ and $P_{\Gamma_M}$, and from them to the growth rate of our graph $\Gamma$.

\begin{proposition}
\label{returndetprop}
Let $\Gamma$ be the graph of some labeled wedge of size $1$ and fix $M$. Let $I$ denote the identity matrix. Then for all $|z| < 2^{\frac{-1}{M}}$ we have:
\[
\frac{P_{\Gamma}(z)}{P_{\Gamma_M}(z)} = \det(I - E(z)).
\]
\end{proposition}

\begin{proof}

Let $c_{n, k}$ be the number of closed paths of length $n$ in $\Gamma$ which visit $S_M$ exactly $k$ times. Let $b_{n, k}$ be the number of these paths which also start in $S_M$. Notice that a path counted in $b_{n, k}$ can be split along its returns to $S_M$ to give $k$ first return paths to $S_M$. Thus:
\[
tr(E(z)^k) = \sum_{n=k}^\infty b_{n, k}z^n.
\]

On the other hand, by definition we know that:
\[
C_\Gamma(n) - C_{\Gamma_M}(n) = \sum_{k=1}^n c_{n, k}.
\]
Consider the number of pairs $(x, j)$ where $x = x_0x_1\dots x_n$ is a path counted by $c_{n, k}$ and $j \in \{0, 1, \dots, n-1 \}$ is a time for which $x_j \in S_M$. By changing the starting vertex of $x$ to $x_j$ we get a pair $(y, j)$ with $y = x_j x_{j+1}\dots x_{j+n}$ (where the indices are $\mod n$) counted by $b_{n, k}$. Moreover, given a pair $(y, j)$ with $y = y_0y_1\dots y_n$ counted by $b_{n, k}$ and $j \in \{0, 1, \dots, n-1 \}$, we can change the starting vertex by $j$ to get a path $x = y_{-j}y_{-j+1}\dots y_{n-j}$ (where again indices are $\mod n$) so that the pair $(x, j)$ is as before. There are $kc_{n, k}$ pairs of the first type and $nb_{n, k}$ pairs of the second type. Thus we get:
\[
kc_{n, k} = n b_{n, k}.
\]
Using this we get that for all $z$ small enough:
\begin{align*}
\sum_{n=1}^\infty \frac{C_\Gamma(n) - C_{\Gamma_M}(n)}{n}z^n &= \sum_{n=1}^\infty \sum_{k=1}^\infty \frac{c_{n, k}}{n}z^n \\
&= \sum_{n=1}^\infty \sum_{k=1}^n \frac{b_{n, k}}{k}z^n \\
&= \sum_{k=1}^\infty \frac{1}{k} \sum_{n=k}^\infty b_{n, k}z^n \\
&= \sum_{k=1}^\infty \frac{tr(E(z)^k)}{k}.
\end{align*}
Thus by definition of spectral determinants we have:
\[
\frac{P_{\Gamma}(z)}{P_{\Gamma_M}(z)} = \exp (-\sum_{k=1}^\infty \frac{tr(E(z)^k)}{k}) = \det(I-E(z)).
\]

\end{proof}

Since $2^{\frac{-1}{M}}$ tends to $1$ as $M$ grows, we have that for any graph $\Gamma$ with inverse growth rate $< 1$ we can always take $M$ large enough to make $2^{\frac{-1}{M}}$ larger than the inverse growth rate. In this case we have that our first return matrix captures the inverse growth rate:

\begin{proposition}
\label{Eprop}
Suppose $r(\Gamma) < 1$. Let $M$ be large enough that $r(\Gamma) < 2^{\frac{-1}{M}}$. Let $\lambda(E_M(t))$ denote the leading eigenvalue of $E_M(t)$. Then for all $0 \leqslant t < r$:
\[
\lambda(E_M(t)) < 1, \quad \lambda(E_M(r)) = 1.
\]
\end{proposition}

\begin{proof}
If for some $t < r(\Gamma)$ we have $\lambda(E_M(t)) \geqslant 1$ then, since the coefficients of $E_M(t)$ vary continuously with $t$, we could find some $t_0 \leqslant t$ with $\lambda(E_M(t_0)) = 1$. This would imply $\det(I-E_M(t_0)) = 0$ and hence, by proposition \ref{returndetprop} that $P_\Gamma(t_0) = 0$ contradicting the fact that $r(\Gamma)$ is the smallest root of $P_\Gamma$.

At $r(\Gamma)$ we know that $P_\Gamma$ is zero. On the other hand, by corollary \ref{nozerocor} we know $P_{\Gamma_M}$ is nonzero. Hence by proposition \ref{returndetprop} we get $\det(I-E_M(r(\Gamma))) = 0$.
\end{proof}

\section{Simplicity of the Entropy Root}
\label{pathcounting}

Here again we fix a size $s = 1$ and consider the graph $\Gamma$ induced by a labeled wedge. The goal of this section is to show that the entropy root of the spectral determinant of graphs induced by unstarred kneading sequences is always simple when $r(\Gamma) < 1$.

\begin{lemma}
\label{periodicitylemma}
Let $w \in \{0, 1, \dots, d-1\}^{\mathbb{N}}$ and let $\Gamma$ be its induced graph. Assume $r(\Gamma) < 1$. Suppose we have a sequence of nonnegative numbers $b_u$ indexed by vertices $u = (i, j)$ of $\Gamma$ such that:
\begin{enumerate}
\item $b_{(1, 2)} = 1$,
\item They satisfy the triangle inequality: $b_{(i, k)} \leqslant b_{(i, j)} + b_{(j, k)}$,
\item $b_u = r(\Gamma) \sum_{u \rightarrow v} b_v$.
\end{enumerate}

Suppose $u_0$ lies on a closed path of $\Gamma$ and $b_{(u_0)} = 0$. Then there is a smallest $p\geqslant 2$ such that $b_{(1, p+1)} = 0$. Moreover, $b_{(i+p, j)} = b_{(i, j)} = b_{(i, j+p)}$ and $b_{(1, j)} = 0$ if and only if $j \equiv 1 \mod p$.
\end{lemma}

\begin{proof}
Since every $b_v$ is nonnegative, the sum rule
\[
b_u = r(\Gamma) \sum_{u \rightarrow v \text{ in } \Gamma} b_v.
\]
implies that if $b_u = 0$ then $b_v = 0$ for any vertex $v$ for which there is a path in $\Gamma$ from $u$ to $v$. In particular, any vertex $v$ reachable from $u_0$ in $\Gamma$ has $b_v = 0$. Since $u_{0}$ lies on a closed path in $\Gamma$, we can follow that path to a separated vertex and follow one of its edges to some vertex $(1, q+1)$. Then using $b_{u_0} = 0$ we obtain $b_{(1, q+1)} = 0$ also. Since such a $q$ exists we can consider $p \geqslant 1$ the smallest natural number such that $b_{(1, p+1)} = 0$. By definition we have $b_{(1, 2)} = 1$ therefore $p \geqslant 2$.

For all $i, j$ we have $b_{(i+1, j+1)} \leqslant \frac{1}{r} b_{(i, j)}$. Indeed, if $(i, j)$ is not separated then the above sum has only one term in it and this is an equality. If the pair is separated then there are two terms in the sum and the triangle inequality gives:
\[
b_{(i+1, j+1)} \leqslant b_{(1, i+1)} + b_{(1, j+1)} = \frac{1}{r}b_{(i, j)}.
\]

Applying this when $b_{(1, p+1)} = 0$ we get $b_{(i, i+p)} = 0$ for $i \geqslant 1$. Using this with the triangle inequality gives:
\[
b_{(i+p, j)} \leqslant b_{(i+p, i)} + b_{(i, j)} = b_{(i, j)}.
\]
But it also gives:
\[
b_{(i, j)} \leqslant b_{(i, i+p)} + b_{(i+p, j)} = b_{(i+p, j)}.
\]
Thus we have a periodicity: $b_{(i+p, j)} = b_{(i, j)}$. Symmetry also gives $b_{(i, j+p)} = b_{(i, j)}$.

Since $p$ is minimal we know $b_{(1, 2)}, b_{(1, 3)}, \dots, b_{(1, p)} > 0$. Thus periodicity gives:
\[
b_{(1, j)} = 0 \iff j \equiv 1 \mod p.
\]

\end{proof}

\begin{lemma}
\label{L1}
Let $w \in \{0, 1, \dots, d-1\}^{\mathbb{N}}$ and let $\Gamma$ be its induced graph. Assume $r(\Gamma) < 1$. Suppose we have a sequence of nonnegative numbers $b_u$ indexed by vertices $u = (i, j)$ of $\Gamma$ such that:
\begin{enumerate}
\item $b_{(1, 2)} = 1$,
\item They satisfy the triangle inequality: $b_{(i, k)} \leqslant b_{(i, j)} + b_{(j, k)}$,
\item $b_u = r(\Gamma) \sum_{u \rightarrow v} b_v$.
\end{enumerate}

Suppose a vertex $u_0$ lies on a closed path in $\Gamma$ and that $b_{u_0} = 0$. Let $p$ be the smallest natural number $\geqslant 2$ such that $b_{(1, p+1)} = 0$ as in lemma \ref{periodicitylemma}. Then:
\[
2r(\Gamma)^p < 1.
\]
\end{lemma}

\begin{proof}
Let $m = \min( b_{(1, j)} \mid 2 \leqslant j \leqslant p) > 0$ and let $j$ be the index for which $m = b_{(1, j)}$.

Consider a path in $\Gamma$, say $v_0v_1\dots v_{p}$ with starting vertex $v_0 = (1, j)$ such that $b_{v_i} > 0$. Each vertex $v_i$ has one or two outgoing edges, and since:
\[
b_u = r(\Gamma) \sum_{u \rightarrow v \text{ in } \Gamma} b_v.
\]
At $v_i$, if there is only one outgoing edge then it is necessarily to $v_{i+1}$ and we have $b_{v_i} = r(\Gamma)b_{v_{i+1}}$. If there are two outgoing edges, say from $v_i$ to $x$ and $y$, then we get:
\[
b_{v_i} = r(\Gamma)(b_x + b_y).
\]
Suppose that for each $v_i$ only one of its outgoing edges $x$ has $b_x > 0$. 

If $v_k = (e, f)$ then if it has only one outgoing edge, it is to $(e+1, f+1)$, both indices have increased by $1$ and hence by $1 \mod p$. If there are two outgoing edges then they are to $(1, e+1)$ and $(1, f+1)$. If one of the targets, say $(1, e+1)$, had $b_{(1, e+1)} = 0$ we have by lemma \ref{periodicitylemma} that $e \equiv 0 \mod p$. Hence, $\mod p$, the indices of $(1, f+1)$ are both $1$ more than the indices of $(e, f)$. Thus as we follow the path from $v_0$ to $v_{p}$, both indices of the vertices are increasing by $1 \mod p$.

By lemma \ref{periodicitylemma}, this would imply:
\[
b_{v_0} = r(\Gamma)^pb_{v_p} = r(\Gamma)^pb_{v_0}
\]
which is impossible as $r(\Gamma) < 1$ and $b_{v_0} > 0$.

Thus for some first vertex $v_i$ with $i < p$ we must have two outgoing edges, say to $x$ and $y$, such that $b_x$ and $b_y$ are both positive.

We now claim that in fact this index $i$ is $< p-1$. Indeed, otherwise if $i=p-1$ then starting at $v_0 = (1, j)$, the indices of $v_{p-1}$ would then be congruent to $(0, j-1) \mod p$. If this vertex has two outgoing edges then one will have a vertex with indices congruent to $1 \mod p$, which by lemma \ref{periodicitylemma} implies the vertex has value $0$, a contradiction.

Let the outgoing edges from $v_{i}$ be to vertices $x$ and $y$. Recall that by definition of $m$ and by lemma \ref{periodicitylemma}, $b_x, b_y \geqslant m$. Applying (3) to our path from $v_0$ to $v_i$ and then at $v_i$ gives:
\[
m = b_{v_0} = r(\Gamma)^i b_{v_i} = r(\Gamma)^{i+1}(b_x + b_y) \geqslant 2r(\Gamma)^{i+1} m.
\]
Then cancelling $m$ from both sides, using the fact that $i+1 < p$ and that $r(\Gamma) < 1$ we get:
\[
2r(\Gamma)^p = (2r(\Gamma)^{i+1})r(\Gamma)^{p-i-1} \leqslant r(\Gamma)^{p-i-1} < 1.
\]
\end{proof}

\begin{lemma}
\label{L2}
Let $w \in \{0, 1, \dots, d-1\}^{\mathbb{N}}$ and let $\Gamma$ be its induced graph. Assume $r(\Gamma) < 1$. Suppose we have a sequence of nonnegative numbers $b_u$ indexed by vertices $u = (i, j)$ of $\Gamma$ such that:
\begin{enumerate}
\item $b_{(1, 2)} = 1$,
\item They satisfy the triangle inequality: $b_{(i, k)} \leqslant b_{(i, j)} + b_{(j, k)}$,
\item $b_u = r(\Gamma) \sum_{u \rightarrow v} b_v$.
\end{enumerate}

Suppose that $u_0$ is some vertex lying on a closed path in $\Gamma$ such that $b_{u_0} = 0$. Then the series: 
\[
B_{u_0}(t) = \sum_{n=0}^\infty B_{\Gamma}(u_0, n)t^n,
\]
converges at $t = (r(\Gamma))$.
\end{lemma}

\begin{proof}
Let $p$ be from lemma \ref{periodicitylemma}. By (3), any vertex $v$ reachable from $u_0$ must also have $b_{v} = 0$.

Suppose $v = (i, j)$ is reachable by a path from $u_0$ and that it is separated. Then there are two outgoing edges to $(1, i+1)$ and $(1, j+1)$. Lemma \ref{periodicitylemma} says that since $b_{(1, i+1)} = b_{(1, j+1)} = 0$, we have $i \equiv j \equiv 0 \mod p$. Hence every separated vertex reachable by a path from $u_0$ has both indices divisible by $p$. But taking a further outgoing edge increases the indices by $1 \mod p$, thus along any path from $u_0$, consecutive separated vertices are at least $p$ edges apart.

In particular, this means for a path of length $N$ starting at $u_0$ there are at most $1 + \lfloor N/p \rfloor$ choices of separated edges. At all other times there is only one outgoing edge. Hence:
\[
B_{\Gamma}(u_0, N) \leqslant 2^{1 + \lfloor N/p \rfloor} \leqslant 2 \cdot 2^{N/p}.
\]
By lemma \ref{L1}, we have $r(\Gamma) 2^{1/p} < 1$ and so:
\[
B_{u_0}(r(\Gamma)) = \sum_{N=0}^\infty B_{\Gamma}(u_0, N)r(\Gamma)^N \leqslant 2\sum_{N=0}^\infty (r(\Gamma)2^{1/p})^N < \infty.
\]
\end{proof}

\begin{proposition}
\label{4prop}
Let $w \in \{0, 1, \dots, d-1\}^{\mathbb{N}}$ and let $\Gamma$ be its induced graph. Assume $r(\Gamma) < 1$. Then $r(\Gamma)$ is a simple root of the spectral determinant $P_{\Gamma}$.
\end{proposition}

\begin{proof}
Let $r = r(\Gamma)$. Choose $M$ large enough that $r< 2^{-1/M}$ and let $E(z)$ be the first return matrix for $S_M$. Then proposition \ref{Eprop} gives $\lambda(E(r)) = 1$.

Consider the graph with vertices $S_M$ and edges between $u$ and $v$ whenever their entry in $E(r)$ is positive. We claim this graph has a single strongly connected component with this growth rate.

Suppose there were two components $C_1$ and $C_2$. Then one of them, say $C_1$ can not be reached from the other. Consider $U$ the set of vertices in $S_M$ which can reach $C_1$. Then $\lambda(E(r) \mid_U) = 1$. Then we can take a right Perron eigenvector $x$ for this restriction. It can be extended to an eigenvector for $E(r)$ by setting it to $0$ when the entry does not belong to $U$. So in particular, $x$ is $0$ for any vertex in $C_2$.

For $u \in S_M$ define $b_u = x_u$. For $u$ outside, define:
\[
b_u = \sum_{\gamma} r^n b_v,
\]
where the sum is taken over all paths $\gamma$ in $\Gamma$ of length $n$ which start at $u$, end at $v \in S_M$ and never visit $S_M$ in between. By lemma \ref{countlemma1} there are at most $4 \cdot 2^{(n-1)/M}$ such paths of length $n$. As there are finitely many $v \in S_M$ and $r \cdot 2^{1/M} < 1$ the sum converges and is bounded uniformly in $u$.

Note that by construction, if $u \notin S_M$ then we have:
\[
b_u = r \sum_{u \rightarrow v} b_v.
\]
When $u \in S_M$ and $v \notin S_M$, since $b_v$ counts first returns to $S_M$, we have that any path counted in $b_v$ can be turned into a first return from $u$ back to $S_M$. Thus we have:
\begin{align*}
b_u &= (E(r) x)_u \\
&= \sum_{\substack{u \rightarrow v \\ v \in S_M}} r x_v
 + \sum_{\substack{u \rightarrow  v \\ v \notin S_M}} r b_v \\
&= r\sum_{u\to v} b_v.
\end{align*}

Also we can define $b_{(i, i)} = 0$. We now check that the $b_u$ satisfy the triangle inequality. Let $\Delta(i, j, k) = b_{(i, k)} - b_{(i, j)} - b_{(j, k)}$ and consider $D = \sup_{i, j, k} \Delta(i, j, k)$. Since the $b_u$ are bounded, we have that $D < \infty$, and taking $i = j$ gives $D \geqslant 0$.

If $w_i = w_j \neq w_k$ then:
\begin{align*}
\frac{\Delta(i, j, k)}{r} &= b_{(1,i+1)}+b_{(1,k+1)}
-b_{(i+1,j+1)}
-b_{(1,j+1)}-b_{(1,k+1)}\\
&=\Delta(i+1,j+1,1)\\
&\leq D.
\end{align*}
The same holds when $w_j = w_k \neq w_i$. In the other cases we have:
\[
\frac{\Delta(i,j,k)}r=
\begin{cases}
\Delta(i+1,j+1,k+1),&w_i=w_j=w_k,\\
\Delta(i+1,1,k+1)-2b_{(1,j+1)},&w_i=w_k\neq w_j,\\
-2b_{(1,j+1)},&w_i,w_j,w_k\text{ all distinct}.
\end{cases}
\]
Each of these are $\leqslant D$. Thus we have $0 \leqslant D \leqslant rD$. As $r < 1$, this means $D = 0$. Hence $\Delta(i, j, k) \leqslant 0$ for all $i, j, k$ and so the $b_u$ satisfy the triangle inequality.

If $b_{(1, 2)} = 0$ then repeatedly applying the triangle inequality and our sum identity would give $b_{(i, j)} = 0$ for all $i, j$. Indeed, if a vertex is non separated then $b_{(i+1, j+1)} = b_{(i, j)} / r$ while if it is separated $b_{(i+1, j+1)} \leqslant b_{(i, j)} / r$. Starting at $(1, 2)$ we can apply this $k-1$ times to get $b_{(k, k+1)} \leqslant r^{1-k}b_{(1, 2)}$. Then the triangle inequality gives that for all $(i, j)$:
\[
b_{(i, j)} \leqslant b_{(i, i+1)} + b_{(i+1, i+2)} + \dots + b_{(j-1, j)} \leqslant b_{(1, 2)} \sum_{k=i}^{j-1}r^{1-k}.
\]
Thus $b_{(1, 2)} > 0$. Rescaling $b_u$ if needed, we can normalize and assume $b_{(1, 2)} = 1$.

Consider the other strongly connected component $C_2$ with inverse growth rate $r$. Then let $F = (E_{uv}(r))_{u, v \in C_2}$. We have $\lambda(F) = 1$. Choose $u_0 \in C_2$. 
Since $C_1$ is not reachable from $C_2$ we have $u_0 \notin U$. As $F$ is irreducible and has leading eigenvalue $1$, we can find a closed path passing through $u_0$ in $\Gamma$. Hence $b_{u_0} = x_{u_0} = 0$. Lemma \ref{L2} then gives:
\[
B_{u_0}(r) < \infty.
\]
By definition we have:
\[
\sum_{j=1}^\infty \sum_{v \in C_2} (F^j)_{u_0v} \leqslant B_{u_0}(r) < \infty.
\]
Let $y$ be a positive right Perron eigenvector of $F$ and normalize so its largest entry is $1$. As $F$ is irreducible, we have all the entries of $y$ are strictly positive. Then:
\[
\sum_{v \in C_2} (F^j)_{u_0v} \geqslant \sum_{v \in C_2} (F^j)_{u_0v} \cdot y_v = (F^jy)_{u_0} = y_{u_0} > 0.
\]
Summing over $j$ gives that $\sum_{j=1}^\infty \sum_{v \in C_2} (F^j)_{u_0v}$ diverges, a contradiction.

Thus our assumption was wrong and no second strongly connected component $C_2$ can exist.

This gives that $1$ is an algebraically simple eigenvalue of $E(r)$. Thus near $r$ we can extend $\lambda(r)$ holomorphically to get:
\[
\det(I-E(z)) = (1 - \lambda(z))g(z),
\]
for some holomorphism $g(z)$ with $g(r) \neq 0$.  For $t < r$ sufficiently close, $\lambda(t)$ is the Perron eigenvalue of $E(t)$. Each first return series has:
\[
E_{uv}(t) = \sum_{n=1}^\infty e_{uv}(n)r^n \big(\frac{t}{r}\big)^n\leqslant \frac{t}{r} \sum_{n=1}^\infty e_{uv}(n)r^n =  \frac{t}{r} E_{uv}(r).
\]
Thus we have:
\[
\lambda(t) \leqslant \frac{t}{r}\lambda(r) = \frac{t}{r}.
\]
So in particular,
\[
\lambda'(r) = \lim_{t \rightarrow r^-} \frac{1-\lambda(t)}{r-t} \geqslant \lim_{t \rightarrow r^-} \frac{1-t/r}{r-t}  = \frac{1}{r} > 0.
\]
As $g(r) \neq 0$ this implies that $z = r$ is a simple root of $\det(I-E(z))$. Since proposition \ref{returndetprop} gives:
\[
\frac{P_{\Gamma}(z)}{P_{\Gamma_M}(z)} = \det(I - E(z)),
\]
and corollary \ref{nozerocor} gives $P_{\Gamma_M}(r) \neq 0$ we have that $r$ is a simple root of $P_{\Gamma}(z)$ also.
\end{proof}

\section{Entropy Gaps}
\label{entropygaps}

The goal of this section is to develop estimates in the difference of entropies. Specifically, for a finite word $v = v_0v_1\dots v_n \in \{0, 1, \dots, d-1\}^{n+1}$, we can let $v^\infty$ denote the infinite sequence of $v$ repeated. Then for a finite word $v$ and a single letter $\tau \in \{0, 1, \dots, d-1 \}$ we wish to compare the entropies of $v^\infty$ and $v\tau^\infty$. 

A useful technical tool in this section will be the internal address. Recall:

\begin{definition}
Given $d \geqslant 2$ and a kneading sequence, define:
\[
\rho_w(k) = \min(n > k \mid w_n \neq w_{n-k}),
\]
with $\rho_w(k) = \infty$ if the set is empty. Then the $\rho_w$ orbit of $k$ is:
\[
k, \rho_w(k), \rho^2_w(k), \dots 
\]
Where we terminate the sequence if the next term will be $\infty$. The internal address of $w$ is the $\rho_w$ orbit of $1$.
\end{definition}

Internal addresses are commonly used when discussing Hubbard trees. As our graph approach is meant to replace the Hubbard tree it should come as no surprise that they are useful here. In fact for our graph they have a very nice interpretation. Recall that at any vertex $(a, b)$ with $a < b$ there are either one or two outgoing edges. In the case of two edges, the one to $(1, a+1)$ is the backwards edge and the other to $(1, b+1)$ is forwards. Thus for any unstarred kneading sequence $w$, starting at $(1, 2)$ there is a unique infinite path in $\Gamma_w$ which always follows forward edges. Then on this path, a vertex $(a, b)$ with $a < b$ is separated if and only if $b$ is in the internal address of $w$. More generally, if we start at a vertex $(1, k+1)$ and take this unique path then a separated vertex $(a, b)$ with $a < b$ occurs on the path if and only if $b$ is in the $\rho_w$ orbit of $k$.

To see why consider the path at $(1, 2)$. By definition, $2 = \rho_w(1)$ if and only if $w_{2} \neq w_{2-1} = w_1$ which is true if and only if $(1, 2)$ is separated. Thus if $2 = \rho_w(1)$ then there is a forward edge to $(1, 3)$. If $2 \neq \rho_w(1)$ then the vertex is not separated and there is only an edge to $(2, 3)$. $3 = \rho_w(1)$ if and only if $w_3 \neq w_{3-1} = w_2$ which is true if and only if $(2, 3)$ is separated. Thus if $3 = \rho_w(1)$ then $(2, 3)$ is separated and there is a forward edge to $(1, 4)$. Otherwise we have one edge to $(3, 4)$. Repeating we see that the separated vertices are exactly the orbit of our starting place under $\rho_w$.

\begin{lemma}
\label{egapl1}
Let $d, p \geqslant 2$ and let $w$ be a periodic kneading sequence with period $p$. Suppose further that $p$ belongs to the internal address of $w$. For every $1 \leqslant k < p$, let:
\[
q(k) = \min(\rho_w(k) - k, p - k).
\]
Then $p$ belongs to the $\rho_w$ orbit of $q(k)$.
\end{lemma}

\begin{proof}
Let $w = (w_1w_2\dots w_p)^\infty$. First suppose that $w_p = \star$. Then it follows $w_j \neq w_p$ for all $j < p$. Thus by definition we have that for all $1 \leqslant j < p$:
\[
j < \rho_w(j) \leqslant p.
\]
In particular, this means $q(k) = \rho_w(k) - k$.

Since the $\rho_w$ orbit of any point is strictly increasing, the above inequality forces that the $\rho_w$ orbit of $j$ contains $p$ for any $1 \leqslant j < p$. As $q(k) < p$ we are done.

Now suppose $w_p \neq \star$. Since $p$ is in the internal address of $w$, \cite{Ka} corollary 3.1.5 implies $w$ has exact period $p$. Thus for any $1 \leqslant k < p$, we have $w_{n-k} \neq w_n$ for some $n > k$, hence $\rho_w(k) < \infty$.

By \cite{Ka} lemma 2.3.14, the $\rho_w$ orbit of $\rho_w(k) - k$ intersects the internal address of $w$. Since $\rho_w(p) = \infty$ and we have assumed the internal address contains $p$, it follows that the address ends at $p$. 

If $\rho_w(k) \leqslant p$ then $q(k) = \rho_w(k) - k$ and the lemma is proven.

If $\rho_w(k) > p$ then $q(k) = p-k$. We claim $\rho_w(q(k)) = \rho_w(k) - k$. For $j \geqslant 1$, periodicity implies $w_j = w_{p+j}$. Thus:
\[
w_{q(k)+j} \neq w_j \iff w_{p+j-k} \neq w_{p+j}.
\]
By definition of $\rho_w(k)$, the symbols on the right agree when $p+j < \rho_w(k)$ and differ when $p+j = \rho_w(k)$. As $\rho_w(k) > p$ the first positive $j$ for which they differ is therefore $j = \rho_w(k) - p$. Then by definition:
\[
\rho_w(q(k)) = q(k) + \rho_w(k) - p = (p-k) + (\rho_w(k) - p) = \rho_w(k) - k.
\]
This implies the orbit of $q(k)$ reaches $\rho_w(k) - k$, so since the orbit of $\rho_w(k) - k$ reaches $p$, it follows that the orbit of $q(k)$ reaches $p$ also.

\end{proof}

\begin{lemma}
\label{egapl2}
Let $d, p \geqslant 2$ and let $w = (w_1w_2\dots w_p)^\infty$ be an unstarred periodic or a $\star$ periodic kneading sequence. Suppose that $p$ belongs to the internal address of $w$. Let $G_w$ be its finite graph and $A_{G_w}$ its adjacency matrix. Suppose a path of length $p-1$ in $G_w$ ends at $(i, p)$.  Then $w_i \neq w_p$.

Moreover, if $x = (x_{(i, j)})$ is a nonnegative left eigenvector of $A_{G_w}$ with nonzero eigenvalue then $x_{(i, p)} = 0$ whenever $w_i = w_p$.
\end{lemma}

\begin{proof}
Let $v_1\dots v_p$ be a path in $G_w$ of length $p-1$ with $v_p = (i, p)$. Notice that at any vertex, if its target has index $t+1$ then the vertex itself has an index $t$. Thus we can write $v_t = (t, j_t)$ for some choice of $j_t$.

By definition, we have $j_{t+1} = j_t + 1$ when $w_{j_t} = w_t$ and $j_t < p$. Otherwise, $j_{t+1} = 1$.

Let $j_1 = k+1$. Then following our path we increase to $j_t = k+t$ until either $j_t = p$ or $w_{j_t} \neq w_t$. The first of these events occurs at:
\[
t= q(k) = \min(\rho_w(k) - k, p - k).
\]
In particular, we have $j_{q(k) + 1} = 1$.

Suppose for some $1 \leqslant a < p$ we have $j_{a+1} = 1$ and $\rho_w(a) \leqslant p$. By definition of $\rho_w$, we have that for any $a < t < \rho_w(a)$, $w_t = w_{t-a}$. Thus, for any $a < t \leqslant \rho_w(a)$, we have $j_t = t-a$. Moreover, if $\rho_w(a) < p$ then $j_{\rho_w(a) + 1} = 1$.

Starting at $a = q(k)$, we can follow its $\rho_w$ orbit to get successive indices $a < p$ with $j_{a+1} = 1$. Note by lemma \ref{egapl1}, $p$ belongs to this orbit. Let $g$ be the last value of this orbit smaller than $p$. Note that by definition, $\rho_w(g) = p$. Then $i = j_p = p-g$. In particular, we have:
\[
w_p \neq w_{p-g} = w_i
\]
which proves the first claim.

For the second claim suppose for some pair $(i, p)$ we have $w_i = w_p$. Then by the above no path of length $p-1$ ends at $(i, p)$. Thus the corresponding column of $A_{G_w}^{p-1}$ is all zeros implying any left eigenvector for a nonzero eigenvalue has a corresponding entry of $0$.
\end{proof}

\begin{lemma}
\label{egapl3}
Let $d, p \geqslant 2$ and let $w = (w_1w_2\dots w_p)^\infty$ be an unstarred periodic or a $\star$ periodic kneading sequence. Let $G_w$ be its finite graph and $A_{G_w}$ its adjacency matrix. Suppose that $A_{G_w}$ has leading eigenvalue $\lambda > 1$ and let $x = (x_{(i, j)})$ be a nonnegative left eigenvector associated to $\lambda$, normalized by:
\[
\sum_{1 \leqslant i < j \leqslant p} x_{(i, j)} = 1.
\]
Then:
\[
\sum_{i=1}^{p-1} x_{(i, p)} = \frac{2(1-\lambda^{-1})\lambda^{1-p}}{1-\lambda^{-p}}.
\]
\end{lemma}

\begin{proof}
For each $j$ let $c_j = \sum_{i \neq j} x_{(i, j)}$. For $2 \leqslant j \leqslant p$, if the target of an edge has index $j$ then its source has index $j-1$. Moreover, any vertex with index $j-1$ has exactly one outgoing edge to a vertex with index $j$. Then since $xA_{G_w} = \lambda x$, we have:
\begin{align*}
\lambda c_j &= \sum_{i \neq j} \lambda x_{(i, j)} \\
&= \sum_{\text{vertices } u} x_u \sum_{i \neq j} (A_{G_w})_{u, (i, j)} \\
&= \sum_{i \neq j-1} x_{(i, j-1)} \\
&= c_{j-1}.
\end{align*}

So in particular, we have that for all $1 \leqslant j \leqslant p$:
\[
c_j = \lambda^{1-j} c_1.
\]
Each $x_{(i, j)}$ occurs exactly twice in the sum $\sum_{j=1}^p c_j$ and so:
\[
2 = \sum_{j=1}^p c_j = \sum_{j=1}^p \lambda^{1-j}c_1 = c_1\frac{1-\lambda^{-p}}{1-\lambda^{-1}}.
\]
Then since $c_p = \lambda^{1-p}c_1$, we obtain:
\[
c_p = \frac{2(1-\lambda^{-1})\lambda^{1-p}}{1-\lambda^{-p}}
\]
as wanted.

\end{proof}

\begin{lemma}
\label{finitetrianglelemma}
Fix $d \geqslant 2$ and let $w$ be an eventually periodic unstarred kneading sequence or a $\star$ periodic kneading sequence. Let $G_w$ be its finite graph. Then for all indices $i, j, k$ and any $n \geqslant 0$,
\[
B_{G_w}((i, k), n) \leqslant B_{G_w}((i, j), n) + B_{G_w}((k, j), n).
\]
\end{lemma}
\begin{proof}
Let $G = G_w$. Fix some $i, j, k$. Let $p$ and $q$ be the period and preperiod for $G_w$. Any index $p+q+1$ should be understood to refer to index $q+1$. We proceed by induction on $n$. The base case $n=0$ is clear as the counts are either $0$ or $1$ depending on if the indices are diagonal or not.

Suppose the claim holds for some $n$. Recall that:
\[
	B_G((i, j), n+1) =
	\begin{cases}
	B_G((i+1, j+1), n) & \text{ if } w_i = w_j, \\
	B_{G}((1, i+1), n) + B_G((1, j+1), n) & \text{ if } w_i \neq w_j.
	\end{cases}
	\tag{*}\label{rec}
\]
If $w_i = w_j = w_k$ then the induction hypothesis and \ref{rec} give the claim.

If $w_i = w_j \neq w_k$ then the induction hypothesis gives:
\[
B_G((1,i+1),n)
\leqslant
B_G((i+1,j+1),n)
+B_G((1,j+1),n).
\]
Adding $B_G((1, k+1), n)$ to both sides and applying \ref{rec} gives the claim.

If $w_i = w_k \neq w_j$ then the induction hypothesis gives:
\[
B_G((i+1,k+1),n)
\leqslant
B_G((1,i+1),n)+B_G((1,k+1),n).
\]
The left hand side of this is exactly in our claimed inequality. Using \ref{rec} on the right hand side of our desired inequality we see it is the same as this right hand side but with two additional nonnegative terms, hence the claimed inequality holds.

The $w_j = w_k \neq w_i$ case is the same as $w_i = w_j \neq w_k$.

Finally if all three of $w_i, w_j$, and $w_k$ are pairwise different then by \ref{rec} the right hand side of our inequality is just the left hand side plus $2B_{G}((1, j+1), n)$ which is nonnegative.

Thus in all cases the inequality holds for $n+1$, hence by induction the claim is true for all $n$.
\end{proof}

\begin{proposition}
\label{egapprop}
Let $d, p \geqslant 2$. Let $P = w_1\dots w_p \in \{0, 1, \dots  , d-1\}^p$ and $w = P^\infty$. Suppose further that $p$ belongs to the internal address of $w$. Suppose $h(w) > 0$ and let $r = e^{-h(w)}$. Then for any $\tau \in \{0, 1, \dots, d-1 \}$ we have:
\[
h(P\tau^\infty) - h(w) \geqslant \frac{(1-r)^2r^{p-1}}{2(1-r^p)}.
\]
\end{proposition}

\begin{proof}
Let $u = P\tau^\infty$. Construct $G_u$ with preperiod $p$ and period $1$. By proposition \ref{tiozzo5.12}, the leading eigenvalue of $A_{G_w}$ is $r^{-1} > 1$. Let $x$ be a nonnegative left eigenvector for $r^{-1}$ normalized so that $\sum_{1 \leqslant i < j \leqslant p} x_{(i, j)} = 1$. Define:
\[
c = \sum_{i=1}^{p-1} x_{(i, p)}.
\]
Then lemma \ref{egapl3} implies:
\[
c = \frac{2(1-r)r^{p-1}}{1-r^p}.
\]
Note also that because we normalized $x$, we have $c \leqslant 1$. Consider:
\[
t = r(1 - \frac{(1-r)c}{4}),
\]
and notice $0 < t < r < 1$. We claim that the leading eigenvalue of $A_{G_u}$ is at least $t^{-1}$. Indeed, suppose $\lambda(A_{G_u}) < t^{-1}$. Define the path series:
\[
B_{(i, j)}(z) = \sum_{n=0}^\infty B_{G_u}((i, j), n)z^n.
\]
Then for any $(i, j)$ we have $B_{(i, j)}(t)$ converges. Recall that by the edge definition we have for $i < p$:
\[
B_{(i, p+1)}(t) = \begin{cases}
1 + tB_{(i+1, p+1)}(t), & w_i = \tau, \\
1 + tB_{(1, i+1)}(t) + tB_{(1, p+1)}(t), & w_i \neq \tau.
\end{cases}
\]
Then if we take the inequality from lemma \ref{finitetrianglelemma}, multiply both sides by $t^n$, and sum over $n$ we get:
\[
B_{(1, i+1)}(t) \leqslant B_{(1, p+1)}(t) + B_{(i+1, p+1)}(t).
\]
Combining this with the above, we have that in either case:
\[
B_{(i, p+1)}(t) \leqslant 1 + tB_{(i+1, p+1)}(t) + 2tB_{(1, p+1)}(t).
\]
At $(p, p+1)$ in $G_u$, if $w_p = \tau$ then there are no outgoing edges. If $w_p \neq \tau$ there are two edges, both to $(1, p+1)$. Thus in this case we would have:
\[
B_{(p, p+1)}(t) = 1 + 2tB_{(1, p+1)}(t).
\]
Let $L = \max_{1 \leqslant i \leqslant p} B_{(i, p+1)}(t)$. Then combining the above with our previous inequality yields:
\[
B_{(i, p+1)}(t) \leqslant 1 + tL + 2tB_{(1, p+1)}(t).
\]
Since this holds for all $i$, taking the maximum gives:
\[
L \leqslant 1 + tL + 2tB_{(1, p+1)}(t),
\]
and hence:
\[
L \leqslant \frac{1 + 2tB_{(1, p+1)}(t)}{1-t}.
\]
Using lemma \ref{finitetrianglelemma} and summing over path lengths as before yields for every $1 \leqslant i < j \leqslant p$:
\begin{equation}
\label{egapBbound}
B_{(i, j)}(t) \leqslant B_{(i, p+1)}(t) + B_{(j, p+1)}(t) \leqslant 2L \leqslant \frac{2 + 4tB_{(1, p+1)}(t)}{1-t}.
\end{equation}
Notice that:
\begin{equation}
\label{egappathsum}
\begin{aligned}
\sum_{1 \leqslant i < j \leqslant p} x_{(i, j)} B_{(i, j)}(t) &= \sum_{1 \leqslant i < j \leqslant p} x_{(i, j)} \left(1 + t\sum_{\substack{(i,j)\to(k,l)\\\text{in }G_u}} B_{(k,l)}(t)\right) \\
&= 1 + t\sum_{1 \leqslant i < j \leqslant p} x_{(i, j)} \sum_{\substack{(i,j)\to(k,l)\\\text{in }G_u}} B_{(k,l)}(t).
\end{aligned}
\end{equation}
When $i, j < p$ the outgoing edges at $(i, j)$ in the graphs $G_u$ and $G_w$ agree. At $(i, p)$, $G_w$ has a single edge to $(1, i+1)$. When $w_i \neq w_p$, $G_u$ has this edge along with a second edge to $(1, p+1)$. When $w_i = w_p$, lemma \ref{egapl2} gives $x_{(i, p)} = 0$. Thus:
\begin{align*}
\sum_{1 \leqslant i < j \leqslant p} x_{(i, j)} \sum_{\substack{(i,j)\to(k,l)\\\text{in }G_u}} B_{(k,l)}(t) &=
\sum_{\substack{1 \leqslant i < p \\ w_i \neq w_p}} x_{(i, p)} B_{(1, p+1)}(t) +
\sum_{1 \leqslant i < j \leqslant p} x_{(i, j)} \sum_{\substack{(i,j)\to(k,l)\\\text{in }G_w}} B_{(k,l)}(t) \\
&= cB_{(1, p+1)}(t) + \sum_{1 \leqslant i < j \leqslant p} x_{(i, j)} \sum_{\substack{(i,j)\to(k,l)\\\text{in }G_w}} B_{(k,l)}(t).
\end{align*}
Since $xA_{G_w} = r^{-1}x$, for each fixed vertex $(k, l)$ of $G_w$ we have:
\[
\sum_{\substack{(i, j) \to (k, l) \\ \text{in } G_w}} x_{(i, j)} = r^{-1}x_{(k, l)}.
\]
Thus by reordering the sums we have:
\begin{align*}
\sum_{1 \leqslant i < j \leqslant p} x_{(i, j)} \sum_{\substack{(i, j) \to (k, l) \\ \text{in } G_w}} B_{(k, l)}(t) &=
\sum_{1 \leqslant k < l \leqslant p} B_{(k, l)}(t) \sum_{\substack{(i, j) \to (k, l) \\ \text{in } G_w}} x_{(i, j)} \\
&= \frac{1}{r} \sum_{1 \leqslant k < l \leqslant p} x_{(k, l)} B_{(k, l)}(t).
\end{align*}
Substituting this back into \eqref{egappathsum} gives:
\[
\sum_{1 \leqslant i < j \leqslant p} x_{(i, j)}B_{(i, j)}(t) = 1 + t \left( cB_{(1, p+1)}(t) + \frac{1}{r} \sum_{1 \leqslant i < j \leqslant p} x_{(i, j)} B_{(i, j)}(t) \right).
\]
Rearranging yields:
\begin{equation}
\label{egapbalance}
(1 - \frac{t}{r}) \sum_{1 \leqslant i < j \leqslant p} x_{(i, j)} B_{(i, j)}(t) = 1 + tcB_{(1, p+1)}(t).
\end{equation}
Using \eqref{egapBbound} and \eqref{egapbalance} along with our choice of $t$ gives:
\begin{align*}
1 + tcB_{(1, p+1)}(t) &= (1 - \frac{t}{r}) \sum_{1 \leqslant i < j \leqslant p} x_{(i, j)} B_{(i, j)}(t) \\
&\leqslant \frac{(1-r)c}{4} \cdot \frac{2 + 4tB_{(1, p+1)}(t)}{1-t} \\
&\leqslant \frac{c}{2} + tcB_{(1, p+1)}(t).
\end{align*}
This implies $1 \leqslant \frac{c}{2}$ contradicting $c \leqslant 1$. Thus our assumption that $\lambda(A_{G_u}) < t^{-1}$ was wrong. In particular, we have:
\[
\lambda(A_{G_u}) \geqslant t^{-1} > r^{-1} > 1.
\]
Thus proposition \ref{tiozzo5.12} gives:
\begin{align*}
h(u) - h(w) &= \log(\lambda(A_{G_u})) + \log(r) \\
&\geqslant \log \frac{r}{t} \\
&= - \log(1 - \frac{(1-r)c}{4}) \\
&\geqslant \frac{(1-r)c}{4} \\
&= \frac{(1-r)^2r^{p-1}}{2(1-r^p)}.
\end{align*}

\end{proof}

\section{Symbolic H\"older Estimates}
\label{secsymbol}

The goal of this section is to prove theorem \ref{symbolicthm}. The proof is broken into two parts, one for the lower bound which is quite straightforward, and a second for the upper bound which will require a few technical estimates. Recall that for two sequences $w, w' \in \{0, 1, \dots, d-1\}^{\mathbb{N}}$ we let $N(w, w') = \min\{i \mid w_i \neq w'_i\}$. Then we define the distance between $w$ and $w'$ as:
\[
D(w, w') = d^{-N(w, w')},
\]
and our local H\"older exponent is:
\[
\alpha_{loc}(h,w) = \sup\left\{ \gamma>0 \mid
\begin{array}{l}
\exists C,\epsilon>0, \text{ if } D(w,w')<\epsilon
\text{ then } \\
|h(w)-h(w')|\leqslant C D(w,w')^\gamma
\end{array}
\right\}.
\]

\begin{proposition}
\label{symboliclower}
Fix $d \geqslant 2$, $H > 0$, and $0 < b < H$. Then there exists a constant $C > 0$ and $N_0$ depending on $d, H,$ and $b$, such that if $w, w'$ both satisfy:
\[
h(w), h(w') \geqslant H,
\]
and if $N(w, w') \geqslant N_0$ then:
\[
|h(w) - h(w')| \leqslant Ce^{-bN(w, w')}.
\]
In particular, for every $w$ with $h(w) > 0$ we have:
\[
\alpha_{loc}(h, w) \geqslant \frac{h(w)}{\log d}.
\]
\end{proposition}

\begin{proof}
Fix $H > 0$ and $0 < b < H$. Choose $M$ large enough that $\frac{\log 2}{M} < H-b$, and choose $a$ such that $e^{-H} < a < e^{-b}2^{-1/M}$. 

For two kneading sequences $w, w'$ let $\Gamma$ and $\Gamma'$ be their respective induced graphs and $E, E'$ be their first return matrices. Let $r = r(\Gamma)$ and $r' = r(\Gamma')$. Note since such graphs have at most two outgoing edges we have:
\[
\frac{1}{2} \leqslant r, r' \leqslant e^{-H} < a.
\]
Set $N_0 = M+1$ and suppose $N(w, w') = N \geqslant N_0$. Then the labels of the graphs $\Gamma$ and $\Gamma'$ agree for any vertex whose indices are $< N$. Thus lemma \ref{unifconv} gives a constant $C_0 = C_0(a, M)$ such that:
\[
\sup_{|z| \leqslant a} ||E(z) - E'(z)||_\infty \leqslant C_0(a2^{1/M})^N.
\]

For any $v \in \{0, 1, \dots, d-1\}^{\mathbb{N}}$ with $h(v) \geqslant H$ let:
\[
r_v = e^{-h(v)}, \quad F_v(z) = \det(I - E_v(z)).
\]
Note by lemma \ref{countlemma1} the entries of $E_v(z)$ are bounded uniformly on $|z| \leqslant a$. 

By proposition \ref{4prop}, $r_v$ is a simple zero of $F_v$. Also $F_v(0) = 1$ and $r_v$ is the smallest root. Thus:
\[
F'_v(r_v) < 0.
\]
Since the sequences $v$ with entropy at least $H$ form a compact set, and the inverse growth rate $r_v$ depends continuously on $v$, there is a constant $c>0$ such that:
\[
F'_v(r_v) \leqslant -2c
\]
for all such $v$.

Then there is some $\delta > 0$ such that if $|t - r_v| < \delta$, $t \leqslant e^{-H}$, then
\[
F'_v(t) \leqslant -c.
\]
Increasing $N_0$ if needed we can assume $N(w, w') \geqslant N_0$ implies $|r-r'| < \delta$.

Then the mean value theorem gives:
\[
c|r-r'| \leqslant |F_w(r') - F_w(r)| = |F_w(r') - F_{w'}(r')|.
\]
Since the entries of the first return matrices are uniformly bounded, we can find a constant $C_1$ so that:
\begin{align*}
|F_w(r') - F_{w'}(r')| &= |\det(I - E(r')) - \det(I - E'(r'))| \\
&\leqslant C_1 ||E(r') - E'(r')||_\infty \\
&\leqslant C_1 \sup_{|z| \leqslant a} ||E(z) - E'(z)||_\infty \\
\leqslant C_1C_0(a2^{1/M})^N.
\end{align*}

Since $r, r' \geqslant 1/2$, the mean value theorem gives:
\[
|h(w) - h(w')| \leqslant 2|r-r'| \leqslant \frac{2C_0C_1}{c}(a2^{1/M})^N \leqslant \frac{2C_0C_1}{c} e^{-bN}.
\]

Fix $w$ with $h(w) > 0$ and choose $0 < b < H < h(w)$. By continuity, all $w'$ close enough to $w$ have $h(w') > H$. Thus applying the estimate above along with:
\[
e^{-bN} = D(w, w')^{\frac{b}{\log(d)}},
\]
we have:
\[
\alpha_{loc}(h, w) \geqslant \frac{b}{\log d}.
\]
Taking $b \rightarrow h(w)$ completes the proof.
\end{proof}

\begin{proposition}
\label{symbolicupper}
Fix $d \geqslant 2$. For every $w \in \{0, 1, \dots, d-1 \}^{\mathbb{N}}$ with $h(w) > 0$, we have:
\[
\alpha_{loc}(h, w) \leqslant \frac{h(w)}{\log d}.
\]
\end{proposition}
\begin{proof}
Fix $w = w_1w_2w_3\dots $ with $h(w) > 0$. For each $p \geqslant 2$ let $q_p$ be the largest entry of the internal address of $w$ that is strictly less than $p$. Note by definition $\rho_w(q_p) \geqslant p$. In particular, $w_n = w_{n-q_p}$ for all $q_p < n < p$. Choose $a_p \in \{0, 1, \dots, d-1 \}$ different from $w_{p-q_p}$ and define: $P_p = w_1w_2\dots w_{p-1}a_p$. Then consider the two words:
\[
v^p = P_p^\infty, \quad u^p = P_p0^\infty.
\]
We claim that $p$ belongs to the internal address of $v^p$. Since $v^p$ and $w$ have the same first $p-1$ digits, their internal addresses agree below $p$, so in particular $q_p$ belongs to the internal address of $v^p$ also. In particular, since $\rho_w(q_p) \geqslant p$, we have $\rho_{v^p}(q_p) \geqslant p$. By our choice of $a_p$, the symbols of $v^p$ at index $p$ and $p-q_p$ differ, hence $\rho_{v^p}(q_p) = p$ and so the internal address of $v^p$ contains $p$.

Suppose for a contradiction that the local H\"older estimate at $w$ holds for some $\gamma$ with:
\[
\gamma > \frac{h(w)}{\log d}.
\]
Then for some constant $C > 0$ and all $w'$ sufficiently close to $w$ we have:
\[
|h(w')-h(w)| \leqslant CD(w, w')^\gamma.
\]
In particular, for all $p$ large enough the above holds for both $w' = v^p$ and $w' = u^p$. Note that $D(w, v^p), D(w, u^p) \leqslant d^{-p}$.

Since $h(v^p) \rightarrow h(w) > 0$ we can take $p$ larger and ensure $h(v^p) > 0$.

Let $r = e^{-h(w)}$ and $r_p = e^{-h(v^p)}$. By proposition \ref{egapprop} we have:
\[
h(u^p) - h(v^p) \geqslant \frac{(1-r_p)^2r_p^{p-1}}{2(1-r_p^p)}.
\]
But on the other hand:
\begin{align*}
h(u^p) - h(v^p) &\leqslant |h(u^p) - h(w)| + |h(v^p) - h(w)| \\
&\leqslant 2Cd^{-\gamma p}.
\end{align*}
So in particular:
\[
\frac{(1-r_p)^2r_p^{p-1}}{2(1-r_p^p)} \leqslant 2Cd^{-\gamma p}.
\]
As $r_p \rightarrow r < 1$, $r_p^p \rightarrow 0$, therefore:
\[
\frac{(1-r_p)^2}{2(1-r_p^p)} \rightarrow \frac{(1-r)^2}{2} > 0.
\]
Thus taking $p$-th roots, as $p \rightarrow \infty$ we get:
\[
r \leqslant d^{- \gamma}.
\]
But $\gamma > \frac{h(w)}{\log d}$ gives: 
\[
d^{-\gamma} < e^{-h(w)} = r,
\]
which is a contradiction.
\end{proof}

Finally, combining propositions \ref{symboliclower} and \ref{symbolicupper} we obtain theorem \ref{symbolicthm}

\section{Angular Approximations}\label{anglelemmas}

Here we develop machinery to turn our symbolic dynamical results above into results on $\Pi(d)$. For our angular lower bound we need only lemma \ref{angularboundarylemma} to control the relationship between two angles being close and their kneading sequences agreeing for a long time by restricting how often the start of the kneading sequence can change. 

The remaining three lemmas concern the upper bound and are mostly aimed at showing that near any starting angle $\theta$ we can find arbitrarily close angles whose kneading sequences are as in proposition \ref{egapprop}.

\begin{lemma}
\label{angularboundarylemma}
Fix $d \geqslant 2$ and $N \geqslant 1$ and let:
\[
B_N = \bigcup_{n=1}^N \bigl\{ \theta \in \mathbb{R} / \mathbb{Z} \mid (d^n-1)\theta \in \mathbb{Z} \bigr\}.
\]
If $I \subseteq \mathbb{R} / \mathbb{Z}$ is an interval of length at most $d^{-N}$, then $I \cap B_N$ contains at most $N$ points. Moreover, for any $\theta, \theta'$ in the same connected component of $I \setminus B_N$, their kneading sequences $K(\theta)$ and $K(\theta')$ agree for the first $N$ symbols and none of these symbols is a $\star$.
\end{lemma}

\begin{proof}
First note for any $n \leqslant N$, the points $\theta$ such that $(d^n-1)\theta \in \mathbb{Z}$ are equally spread out with distance:
\[
\frac{1}{d^n - 1} > \frac{1}{d^N}.
\]
So in particular $I$ contains at most one of these points for each $n$. Thus $I \cap B_N$ has at most $N$ points.

The $n$-th symbol of $K(\theta)$ is a $\star$ when $d^{n-1}\theta$ belongs to $M_\theta$. This happens when, for some $j < d$:
\[
d^{n-1}\theta = \frac{\theta + j}{d} \mod 1,
\]
which is the same as:
\[
(d^n-1)\theta \in \mathbb{Z}.
\]
Thus by definition of $B_N$, no element not in $B_N$ has a $\star$ in its kneading sequences first $N$ places.

Consider $J$ a connected component of $I \setminus B_N$. As $\theta$ varies in $J$, its first $N$ orbit points vary continuously on $\mathbb{R} / \mathbb{Z}$. No point meets an endpoint of $M_\theta$ (as this would put a $\star$ in the kneading sequence), thus these first $N$ iterates of $\theta$ have the same label relative to $M_\theta$. 
\end{proof}

\begin{lemma}
\label{angularaddresslemma}
Fix $d \geqslant 2$ and $m \geqslant 1$. Let $\theta = \frac{k}{d^m} \in (0, 1)$ with base $d$ expansion:
\[
\theta = 0.a_1a_2\dots a_m, \quad a_m \neq 0.
\]
Let $w = K(\theta)$. Let $R$ be the length of the longest consecutive block in $a_1\dots a_m$ of all $0$'s or all $d-1$'s. If neither digit occurs take $R = 0$. Then every $p \geqslant m+R+2$ belongs to the internal address of $w$. 
\end{lemma}

\begin{proof}
The orbit of $\theta$ reaches $0$ after exactly $m$ iterations. In particular this means $\theta \neq 0$. Since the region labeled by $d-1$ contains $0$, $w_j = d-1$ for all $j > m$. Notice $d^{m-1}\theta = \frac{a_m}{d} \mod 1$. Since $a_m \neq 0$, we have:
\[
\theta < 1 \leqslant a_m \leqslant d-1 < \theta + d - 1.
\]
Since the region labeled by $d-1$ is:
\[
[0, \frac{\theta}{d}) \cup (\frac{\theta + d - 1}{d}, 1),
\]
we have $w_1, w_m \neq d-1$. Suppose for some $i \leqslant j < m$ we have $w_i, w_{i+1}, \dots, w_{j}$ are all $d-1$. If $d^{i-1}\theta \in [0, \frac{\theta}{d})$ then $d^{i}\theta \in [0, \theta)$ which has empty intersection with $(\frac{\theta + d - 1}{d}, 1)$. Thus for $w_{i+1} = d-1$ we would have $d^{i}\theta \in [0, \frac{\theta}{d})$ also. Similarly, if $d^{i-1}\theta \in (\frac{\theta+d-1}{d}, 1)$ then $d^{i}\theta \in (\theta, 1)$ which has empty intersection with $[0, \frac{\theta}{d})$. Thus for $w_{i+1} = d-1$ we would have $d^{i}\theta \in (\frac{\theta + d - 1}{d}, 1)$.

Repeating this, we have that for each index $i, i+1, \dots, j$ either all the corresponding iterates of $\theta$ must belong to $[0, \frac{\theta}{d})$ or all must belong to $(\frac{\theta + d - 1}{d}, 1)$. This implies that the corresponding digits $a_i, a_{i+1}, \dots, a_{j}$ are either all $0$ or all $d-1$. Thus we have $j-i +1 \leqslant R$.

For $j \geqslant m$, since the tail of $w$ is constantly $d-1$ and $w_1 \neq d-1$, we have:
\[
\rho_w(j) = j+1.
\]
For $1 \leqslant j < m$, since $w_{m+j} = d-1 \neq w_m$ we have:
\[
\rho_w(j) \leqslant m + j.
\]
If $\rho_w(j) > m$ then by definition, for any $m+1-j \leqslant k < \rho_w(j) - j$ we have $w_k = d-1$. There are $\rho_w(j) - m - 1$ consecutive symbols here, thus $\rho_w(j)-m-1 \leqslant R$. Put another way:
\[
\rho_w(j) \leqslant m + R + 1.
\]
Since this holds for every $j \leqslant m$, and when $j \geqslant m$, $\rho_w(j)$ increases by $1$, we are done.
\end{proof}

\begin{lemma}
\label{angularcompanionlemma}
Fix $d \geqslant 2$ and $m \geqslant 1$. Let $\theta = \frac{k}{d^m} \in (0, 1)$ with base $d$ expansion:
\[
\theta = 0.a_1a_2\dots a_m, \quad a_m \neq 0.
\]
Let $w = K(\theta)$. Let $R$ be the length of the longest consecutive block in $a_1\dots a_m$ of all $0$'s or all $d-1$'s. If neither digit occurs take $R = 0$. For every $p \geqslant m+R+2$, define:
\[
\phi = \frac{\theta}{1-d^{-p}}.
\]
Then $\phi$ has period $p$ under multiplication by $d$ and:
\[
K(\phi) = (w_1w_2\dots w_{p-1}\star)^\infty.
\]
\end{lemma}

\begin{proof}
First note that $\phi \in (0, 1)$. Indeed, since $ \theta = k/d^m \in (0, 1)$ and $k$ is an integer, we have:
\[
0 < \theta \leqslant 1 - \frac{1}{d^m}.
\]
As $p > m$, $d^{-p} < d^{-m}$ so:
\[
0 < \phi = \frac{\theta}{1-d^{-p}} \leqslant \frac{1-d^{-m}}{1-d^{-p}} < 1
\]
Rearranging the definition of $\phi$ gives:
\[
\phi - \theta = \phi d^{-p}.
\]
So we obtain:
\[
0 < \phi - \theta < d^{-p}.
\]

First we show that the kneading sequences $K(\theta)$ and $K(\phi)$ agree for their first $m$ symbols. Fix $n \leqslant m$ and pick $b > \theta$ so that $(d^n-1)b \in \mathbb{Z}$. Then we can write $b = \frac{s}{d^n-1}$. Notice:
\[
b - \theta = \frac{sd^m - k(d^n-1)}{d^m(d^n-1)}.
\]
Since $n \leqslant m$ both $sd^m$ and $kd^n$ are divisible by $d^n$. Therefore the numerator has the same remainder as $k$ does when divided by $d^n$. Using the base $d$ expansion of $\theta$ we know that:
\[
k = \sum_{j=1}^m a_jd^{m-j}.
\]
The terms from $j \leqslant m-n$ are divisible by $d^n$ so the remainder is:
\[
\sum_{j=m-n+1}^m a_jd^{m-j}.
\]
In particular, its base $d$ expansion is $a_{m-n+1}a_{m-n+2}\dots a_m$. Since there are at most $R$ consecutive zeros, this number is at least $d^{n-R-1}$. Thus we have:
\[
sd^m - k(d^n-1) \geqslant d^{n-R-1}.
\]
Since $p \geqslant m + R + 2$, this implies:
\[
\phi - \theta < d^{-p} < d^{-m-R-1} = \frac{d^{n-R-1}}{d^md^n} < \frac{d^{n-R-1}}{d^m(d^n-1)} \leqslant b - \theta.
\]
In particular we have $b > \phi$. Let:
\[
B_m = \bigcup_{n=1}^m \{ x \in \mathbb{R} / \mathbb{Z} \mid (d^n - 1)x \in \mathbb{Z} \}.
\]
Then the above inequality gives $(\theta, \phi] \cap B_m = \emptyset$. Note also that $\theta \notin B_m$ as $d^m\theta = 0 \neq \theta \mod 1$ so the orbit of $\theta$ is not periodic.

Since the interval $[\theta, \phi]$ has length at most $d^{-p} \leqslant d^{-m}$, lemma \ref{angularboundarylemma} implies that $K(\theta)$ and $K(\phi)$ agree in their first $m$ symbols and none are $\star$.

For $n > m$ we have $w_n = d-1$ as $d^m\theta = 0 \mod 1$. Assume that $n \leqslant p$. Then:
\[
d^{n-1}\phi = d^{n-1}\theta + \phi d^{n-1-p} = \phi d^{n-1-p} \mod 1.
\]
Note $0 < \phi d^{n-1-p} \leqslant \phi/d < 1$. For $n < p$ the second inequality is strict so the $n$-th symbol of $K(\phi)$ is $d-1$.  When $n=p$ this gives that $d^{n-1}\phi = \phi/d \mod 1$ which is an endpoint of the leaf of $M_\phi$.  Thus the kneading sequence contains a $\star$ at its $p$-th entry. In particular, we have:
\[
K(\phi) = (w_1w_2\dots w_{p-1}\star)^\infty.
\]
In particular, since the $\star$ first occurs at entry $p$, $\phi$ has period $p$ under multiplication by $d$.
\end{proof}

\begin{lemma}
\label{angularapproximationlemma}
Fix $d \geqslant 2$ and $\theta \in (0, 1)$. There exists a $C > 0$ such that for arbitrarily large integers $p$, there is an angle $\theta' \in (0, 1)$ with base $d$ expansion:
\[
\theta' = 0.b_1b_2\dots b_m, \quad b_m \neq 0,
\]
such that:
\[
m+R+2 \leqslant p, \quad |\theta' - \theta| \leqslant Cd^{-p}.
\]
Here $|\theta'-\theta|$ is taken over $\mathbb{R} / \mathbb{Z}$ and $R$ is the length of the longest consecutive block in $b_1\dots b_m$ of all zeros or all $d-1$'s. If neither digit occurs take $R = 0$.

\end{lemma}

\begin{proof}
If $\theta$ has a finite base $d$ expansion we can take $\theta' = \theta$. Suppose this is not the case. Then $\theta$ has some base $d$ expansion:
\[
\theta = 0.a_1a_2\dots,
\]
which is neither eventually all zeros nor all $d-1$'s.

First suppose that the length of any block of consecutive zeros or $d-1$'s in $a_1a_2\dots$ is bounded by some $B$. Then there are arbitrarily large indices $m$ with $a_m \neq 0$. Take:
\[
\theta' = \sum_{j=1}^m a_jd^{-j},
\]
Note for this choice of $\theta'$, $R \leqslant B$. For $p = m +R + 2$ we have:
\[
0 \leqslant  \theta - \theta' \leqslant \sum_{j=m+1}^\infty (d-1)d^{-j} = d^{-m} \leqslant d^{B+2}d^{-p}.
\]
In particular, we have $|\theta' - \theta| \leqslant d^{B+2}d^{-p}$. Take $C = d^{B+2}$. Since we can find arbitrarily large $m$ we get arbitrarily large $p$ satisfying the lemma statement so we are done.

Finally we consider the case that the lengths of consecutive blocks of zeros or $d-1$'s are unbounded. Consider a maximal block of all zeros or $d-1$'s of length $L$ such that any earlier block has length strictly less than $L$. Note we can do this for arbitrarily large $L$. Assume this block occurs at index $m+1$. Let $Q$ be the length of the longest block of zeros or $d-1$'s in $a_1a_2\dots a_m$. Then by assumption $L \geqslant Q+1$.

If the block $a_{m+1}\dots a_{m+L}$ consists of zeros, take:
\[
\theta' = \sum_{j=1}^m a_jd^{-j}.
\]
Since the block was assumed to be maximal, $a_m \neq 0$. Note $R = Q$. Moreover:
\[
0 \leqslant \theta - \theta' = \sum_{j=m+L+1}^\infty a_jd^{-j} \leqslant \sum_{j=m+L+1}^\infty (d-1)d^{-j} = d^{-(m+L)}.
\]
If the block $a_{m+1}\dots a_{m+L}$ consists of $d-1$'s take:
\[
\theta' = \sum_{j=1}^m a_jd^{-j}+d^{-m}.
\]
As the block is maximal, $a_m < d-1$ so the added $d^{-m}$ changes the final digit from $a_m$ to $a_{m}+1$. In particular, $\theta' \in (0, 1)$ and has a length $m$ base $d$ expansion. Changing the final digit could increase the maximal length of a block consisting of only $d-1$'s by $1$, hence $R \leqslant Q+1$.

Since we assumed $a_{m+1}\dots a_{m+L}$ are all $d-1$ we have:
\begin{align*}
\theta' - \theta &= d^{-m} - \sum_{j=m+1}^{m+L} (d-1)d^{-j}  - \sum_{j=m+L+1}^\infty a_jd^{-j} \\
&= d^{-(m+L)} - \sum_{j=m+L+1}^\infty a_jd^{-j} \\
\end{align*}
Since the sum is between $0$ and $d^{-(m+L)}$ we have:
\[
0 \leqslant \theta' - \theta \leqslant d^{-(m+L)}.
\]
Thus regardless of if the block $a_{m+1}\dots a_{m+L}$ is zeros or $d-1$'s we can construct $\theta'$ with $R \leqslant Q+1$ and $|\theta' - \theta| \leqslant d^{-(m+L)}$. Taking $p = m + R + 2$ and using $L \geqslant Q+1$ gives $p \leqslant m + Q + 3 \leqslant m + L + 2$ and so:
\[
|\theta' - \theta| \leqslant d^{-(m+L)} = d^{p-m-L}d^{-p} \leqslant d^2d^{-p}.
\]
Since we can take $m$ to be arbitrarily large by selecting blocks with larger $L$, we can take $p$ arbitrarily large which completes the proof.
\end{proof}

\section{The Local H\"older Exponent}\label{mainproofsection}

The goal of this section is to prove theorem \ref{thm}. The approach is to combine our already proven theorem \ref{symbolicthm} and use the tools developed in section \ref{anglelemmas} to convert this symbolic estimate into an estimate on $\Pi(d)$. As with the proof of theorem \ref{symbolicthm}, we break our proof into two parts, one for the lower bound and the other for the upper.

\begin{proposition}
\label{angularlower}
Fix $d \geqslant 2$. For any $M \in \Pi(d)$ with $h(M) > 0$, we have:
\[
\alpha_{loc}(h\mid_{\Pi(d)}, M) \geqslant \frac{h(M)}{\log d}.
\]
\end{proposition}

\begin{proof}
Let $M = M_\theta$ and fix some $0 < \gamma < h(M)/\log d$. Pick $b$ and $H$ so that:
\[
\gamma \log d < b < H < h(M).
\]
Then by continuity we can take some interval $U$ around $\theta$ such that for any $\eta \in U$, $h(M_\eta) > H$. If $\eta$ is irrational then it has an unstarred kneading sequence and $h(K(\eta)) = h(M_\eta)$. Thus by proposition \ref{symboliclower} there is some constant $C > 0$ and $N_0$ such that for any irrational angles $\eta, \eta' \in U$, if the kneading sequences of $\eta$ and $\eta'$ agree in their first $N \geqslant N_0$ symbols then:
\[
|h(M_\eta) - h(M_{\eta'})| \leqslant Ce^{-bN}.
\]
Let $\theta' \neq \theta$ and set $\delta = |\theta' - \theta|$ where the absolute value is taken over $\mathbb{R} / \mathbb{Z}$. Taking $\theta'$ close enough to $\theta$ we can choose $N$ such that $N \geqslant N_0$ and:
\[
d^{-(N+1)} < \delta \leqslant d^{-N}.
\]
Note also that for $\delta$ small enough the smallest interval $I$ joining $\theta$ to $\theta'$ on $\mathbb{R} / \mathbb{Z}$  will lie in $U$.

By lemma \ref{angularboundarylemma}, the set:
\[
B_N = \bigcup_{n=1}^N \{\eta \in \mathbb{R} / \mathbb{Z} \mid (d^n-1)\eta \in \mathbb{Z} \}
\]
contains at most $N$ points of $I$. We can list $\theta$, $\theta'$, and the points of $B_N$ in $I$ ordered from $\theta$ to $\theta'$ as:
\[
\xi_0 = \theta, \xi_1, \dots, \xi_q = \theta'.
\]
Note $q \leqslant N+1$. For each $j = 1, \dots, q$, lemma \ref{angularboundarylemma} gives that all angles in the  smallest open interval from $\xi_{j-1}$ to $\xi_j$ have the same first $N$ symbols in their kneading sequences. Thus for any irrational angles $\eta, \eta'$ in this interval, we have:
\[
|h(M_\eta) - h(M_{\eta'})| \leqslant Ce^{-bN}.
\]
Then by continuity we can approach the endpoints of the intervals to get that for all $j = 1, \dots, q$:
\[
|h(M_{\xi_j}) - h(M_{\xi_{j-1}})| \leqslant Ce^{-bN}.
\]
Since $\xi_0 = \theta$ and $\xi_q = \theta'$ we have:
\begin{align*}
|h(M_{\theta'}) - h(M_{\theta})| &= |\sum_{j=1}^q \left( h(M_{\xi_j}) - h(M_{\xi_{j-1}}) \right) | \\
&\leqslant \sum_{j=1}^q |h(M_{\xi_j}) - h(M_{\xi_{j-1}})| \\
&\leqslant C(N+1)e^{-bN}
\end{align*}
Notice:
\[
C(N+1)e^{-bN} = C(N+1)e^{-(b-\gamma \log d)N}d^{-\gamma N}.
\]
As $b > \gamma \log d$, the term $C(N+1)e^{-(b-\gamma \log d)N}$ is bounded. Thus we obtain:
\[
|h(M_{\theta'}) - h(M_{\theta})| \leqslant C' d^{-\gamma N}.
\]
Since we chose $d^{-(N+1)} < |\theta - \theta'|$ we obtain:
\[
|h(M_{\theta'}) - h(M_{\theta})| \leqslant C' d^\gamma |\theta - \theta'|^\gamma.
\]
Since $C'd^\gamma$ is independent of $\theta'$, we have $\alpha_{loc}(h\mid_{\Pi(d)}, M) \geqslant \gamma$ for any $\gamma < h(M)/\log d$ as wanted.
\end{proof}

\begin{proposition}
\label{angularupper}
Fix $d \geqslant 2$. For any $M \in \Pi(d)$ with $h(M) > 0$ we have:
\[
\alpha_{loc}(h\mid_{\Pi(d)}, M) \leqslant \frac{h(M)}{\log d}.
\]
\end{proposition}

\begin{proof}
Let $M = M_\theta$ for $\theta \in (0, 1)$. Suppose for a contradiction that the upper bound is false at $M$. Then we can find a $\gamma > h(M) / \log d$ such that for a constant $C > 0$, any angle $\eta$ close to $\theta$ gives:
\[
|h(M_\eta) - h(M)| \leqslant C|\eta - \theta|^\gamma.
\]
Lemma \ref{angularapproximationlemma} gives a constant $C_0 > 0$ such that, for arbitrarily large $p$, we can pick a $\theta' \in (0, 1)$ with base $d$ expansion:
\[
\theta' = 0.b_1b_2\dots b_m, \quad b_m \neq 0,
\]
with the property that $|\theta' - \theta| \leqslant C_0d^{-p}$ and $m+R+2 \leqslant p$. Here $R$ denotes the length of the largest consecutive block of all zeros or all $d-1$'s in $b_1b_2\dots b_m$.

Let $w = K(\theta')$ and $\phi = \frac{\theta'}{1-d^{-p}}$. Then lemma \ref{angularcompanionlemma} gives $\phi \in (0, 1)$ and:
\[
K(\phi) = (w_1w_2\dots w_{p-1} \star)^\infty.
\]
Note also that $0 < \phi - \theta' < d^{-p}$. Thus:
\[
|\phi - \theta| \leqslant |\phi - \theta'| + |\theta' - \theta| \leqslant (C_0+1)d^{-p}.
\]
Let $P = w_1w_2\dots w_p$ and $v = P^\infty$. Note that since $d^m\theta' = 0 \mod 1$ and we label the sector of the circle containing $0$ by $d-1$, we have $w_j = d-1$ for all $j > m$. As $p > m$, we have $w = P(d-1)^\infty$. By lemma \ref{angularaddresslemma}, $p$ is in the internal address of $w$. As $w$ and $v$ agree for their first $p$ symbols, $p$ is also in the internal address of $v$.

Note the labelled wedges induced by $M_\phi$ and $v$ are periodic of period $p$, and their labels agree at any index $(i, j)$ provided neither $i$ nor $j$ is divisible by $p$. Lemma \ref{starcomparison} gives: 
\[
h(M_\phi) = h(v).
\]
Note by definition we have $h(M_{\theta'}) = h(w)$.

Since $h(v) \rightarrow h(M) > 0$, by continuity we have for all $p$ sufficiently large, $h(v) > 0$ also. Let $r_p = e^{-h(v)}$. Then proposition \ref{egapprop} with $\tau = d-1$ gives:
\[
h(M_{\theta'}) - h(M_\phi) = h(P(d-1)^\infty) - h(P^\infty) \geqslant \frac{(1-r_p)^2r_p^{p-1}}{2(1-r_p^p)}.
\]
But we also know that:
\begin{align*}
h(M_{\theta'}) - h(M_\phi) &\leqslant |h(M_{\theta'}) - h(M)| + |h(M_\phi) - h(M)| \\
&\leqslant C|\theta' - \theta|^\gamma + C|\phi - \theta|^\gamma \\
&\leqslant 2C(C_0+1)^\gamma d^{-\gamma p}.
\end{align*}
Combining these gives:
\[
\frac{(1-r_p)^2r_p^{p-1}}{2(1-r_p^p)} \leqslant 2C(C_0+1)^\gamma d^{-\gamma p}.
\]
Taking $p$-th roots gives:
\[
\left( \frac{(1-r_p)^2}{2(1-r_p^p)} \right)^{1/p} r_p^{1-1/p} \leqslant (2C(C_0+1)^\gamma)^{1/p}d^{-\gamma}.
\]
Let $r = e^{-h(M)}$. As $h(M) > 0$ we have $r < 1$, so since $r_p \rightarrow r$ we have $r_p^p \rightarrow 0$. Thus:
\[
\frac{(1-r_p)^2}{2(1-r_p^p)} \rightarrow \frac{(1-r)^2}{2} > 0.
\]
Thus taking limits we obtain:
\[
r \leqslant d^{-\gamma}.
\]
Since $\gamma > h(M) / \log d$, $d^{-\gamma} < e^{-h(M)} = r$ giving a contradiction.
\end{proof}

Finally, combining propositions \ref{angularlower} and \ref{angularupper} we obtain theorem \ref{thm}.

\section{A Multicritical Counterexample}\label{counterexmaplesection}

Finally we show that the general conjecture $\alpha_{loc}(h, M) = h(M) / \log d$ is false when $M$ is multicritical. While much of the machinery developed in earlier sections is true for multicritical primitive majors, the fact that the entropy root of the spectral determinant is simple from proposition \ref{4prop} fails to generalize. In fact conceptually it is simple to find counterexamples to this estimate. Indeed, consider the family of primitive majors in $PM(3)$ which this entire section will focus on:
\[
M_N = \big\{ \{1/9 - \delta_N, 4/9 -\delta_N\}, \{31/54, 49/54\} \big\}, \quad \delta_N = \frac{1}{6 \cdot 3^N},
\]
and their limit:
\[
M = \big\{ \{1/9, 4/9 \}, \{ 31/54, 49/54\} \big\}.
\]
These majors are constructed to have the following property. For $M$, the orbits of the two leaves are:
\[
1/3 \mapsto 0, \quad 13/18 \mapsto 1/6 \mapsto 1/2.
\]
Thus each leaf eventually maps to a distinct fixed point of multiplication by $3$. Informally, this produces a higher multiplicity for the growth rates in the graph of $M$. On the other hand, for $M_N$, the addition of $\delta_N$ makes both leaves map to the fixed point $1/2$ giving vastly different behaviour in the orbit of the leaves for majors arbitrarily close to $M$.

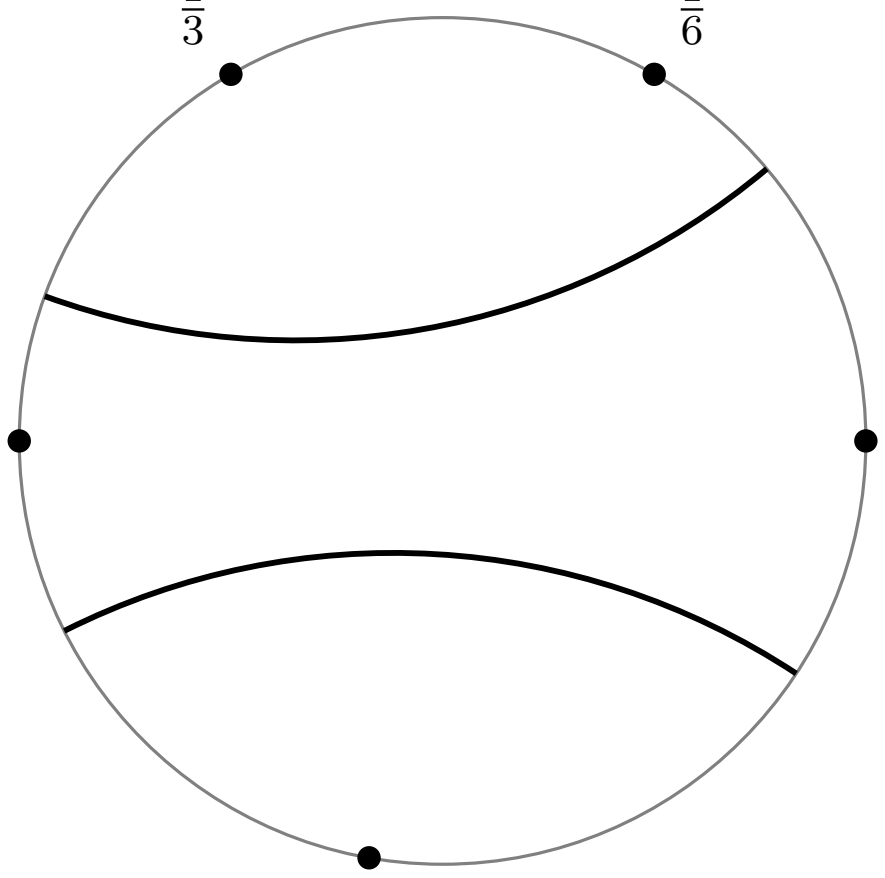
\begin{figure}[htbp]
\centering
\begin{tikzpicture}[font=\normalsize]
  \draw[black!50,line width=.6pt] (0,0) circle[radius=2.8cm];

  \draw[line width=1.1pt] (40:2.8)
    arc[start angle=-50,end angle=-110,radius=4.849742cm];
  \draw[line width=1.1pt]
    ({620/3}:2.8)
    arc[start angle={350/3},end angle={170/3},radius=4.849742cm];

  \foreach \angle/\value in {
    0/{0},
    60/{\frac16},
    120/{\frac13},
    180/{\frac12},
    260/{\frac{13}{18}}
  } {
    \fill (\angle:2.8) circle[radius=2.2pt];
    \node at (\angle:3.3) {$\value$};
  }
\end{tikzpicture}
\caption{The primitive major $M$ with leaves $\{1/9, 4/9\}$ and $\{31/54, 49/54\}$, and its postcritical set.}
\label{majorfig}
\end{figure}

First we recall a finite graph model for core entropy of primitive majors whose leaves have rational endpoints (see \cite{GT}). Given such a major $M = \{ L_1, L_2 \}$ we can let $P$ be the forward orbit of the endpoints of the leaves under $T(x) = 3x \mod 1$. If the endpoints are rational numbers then this set is finite. By abuse of notation we let $T(L_j)$ denote the image of any endpoint of the leaf $L_j$ under $T$. By definition all endpoints map to the same point so this is well defined. Then we define a graph $G_M$ with vertex set all unordered pairs of distinct elements of $P$. 

Given a vertex $(x, y)$ suppose that the arc from $x$ to $y$ crosses no leaves. Then there is one edge to $(T(x), T(y))$. If the vertex passes through only one leaf $L$ then we have two edges from $(x, y)$. One to $(T(L), T(x))$ and the other to $(T(L), T(y))$. Finally if we pass through both leaves in the order $L_{a_1}$ then $L_{a_2}$ when traveling from $x$ to $y$ then we have three edges from $(x, y)$. One to $(T(L_{a_1}), T(x))$, one to $(T(L_{a_1}), T(L_{a_2}))$, and finally to $(T(L_{a_2}), T(y))$. In our edge rule we discard any targets which would be diagonal. Also we allow multiple edges with the same source and target.

Then if $A_M$ denotes the adjacency matrix of this graph, and if $h(M) > 0$, we have $h(M) =  \log \lambda(A_M)$.

\begin{figure}[htbp]
\centering
\begin{tikzpicture}[
  >={Stealth[length=2mm,width=1.5mm]},
  vertex/.style={ellipse,draw=black!75,fill=white,
    minimum width=17mm,minimum height=10mm,inner xsep=2mm,
    inner ysep=1mm,font=\normalsize},
  every edge/.style={draw=black!80,->,semithick},
  component/.style={draw=black!45,dashed,rounded corners=4pt,
    fill=black!2,inner sep=4mm},
  x=1cm,y=1cm]
  \node[vertex] (B) at (0,1)
    {$\bigl(\frac16,\frac{13}{18}\bigr)$};
  \node[vertex] (C) at (0,3.4)
    {$\bigl(\frac13,\frac{13}{18}\bigr)$};
  \node[vertex] (D) at (4.4,1)
    {$\bigl(\frac13,\frac12\bigr)$};
  \node[vertex] (A) at (4.4,3.4)
    {$\bigl(0,\frac13\bigr)$};
  \node[vertex] (E) at (-3.3,3.4)
    {$\bigl(\frac12,\frac{13}{18}\bigr)$};
  \node[vertex] (F) at (-3.3,1)
    {$\bigl(0,\frac{13}{18}\bigr)$};
  \node[vertex] (P) at (7.4,2.2)
    {$\bigl(0,\frac16\bigr)$};
  \node[vertex] (Q) at (2.3,-1.5)
    {$\bigl(\frac16,\frac12\bigr)$};
  \node[vertex] (H) at (-0.4,-1.5)
    {$\bigl(0,\frac12\bigr)$};
  \node[vertex] (J) at (-3.3,-1.5)
    {$\bigl(\frac16,\frac13\bigr)$};
  \coordinate (BCtop) at (0,4.6);
  \coordinate (BCbottom) at (0,-0.2);
  \coordinate (Atop) at (4.4,4.6);
  \begin{scope}[on background layer]
    \node[component,fit=(B)(C)(BCtop)(BCbottom)] {};
    \node[component,fit=(A)(Atop)] {};
  \end{scope}
  \node[font=\small] at (0,5.3) {growth $2$};
  \node[font=\small] at (4.4,5.3) {growth $2$};
  \path
    (B) edge[loop below,min distance=8mm] (B)
    (C) edge[loop above,min distance=8mm] (C)
    (B) edge[bend left=20] (C)
    (C) edge[bend left=20] (B)
    (B) edge (D)
    (C) edge (A)
    (D) edge (A)
    (D) edge[out=-35,in=-85,loop,looseness=5,min distance=8mm] (D)
    (A) edge[loop above,min distance=8mm]
      node[above=4pt,pos=.5,font=\small,inner sep=1pt] {$2$} (A)
    (E) edge[loop above,min distance=8mm] (E)
    (F) edge[loop left,min distance=8mm] (F)
    (H) edge[loop below,min distance=8mm] (H)
    (E) edge (B)
    (F) edge (B)
    (P) edge (A)
    (P) edge (D)
    (Q) edge node[below right=5pt,pos=.48,font=\small,
      inner sep=1pt] {$2$} (D)
    (J) edge (H);
\end{tikzpicture}
\caption{The entire graph $G_M$. Dashed boxes are the two
components with growth rate $2$. The arrows marked $2$ represent two edges with the same source and target.}
\label{graphfig}
\end{figure}
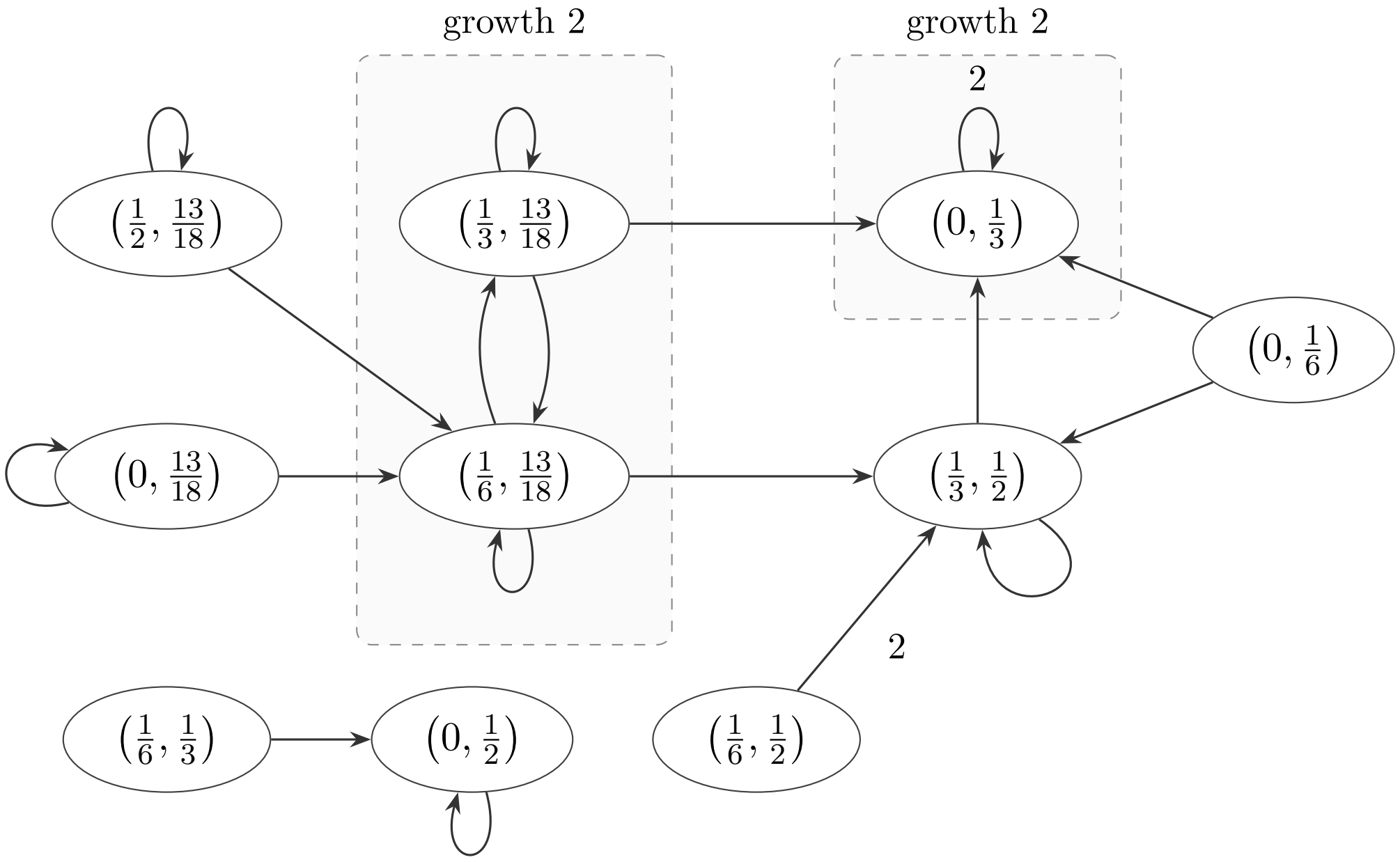

For instance, in our example $M$ above we have $P = \{0, 1/6, 1/3, 1/2, 13/18 \}$. See figure \ref{graphfig} for the finite graph. Based on the figure it is clear that $h(M) = \log 2$ as there are two strongly connected components of the graph, each with growth rate $2$, and no other strongly connected component has growth rate above $2$.

\begin{lemma}
\label{cubicmetriclemma}
For each $N \in \mathbb{N}$, we have:
\[
md(M_N, M) \leqslant 2\delta_N = 3^{-(N+1)}.
\]
\end{lemma}

\begin{proof}
Note that $M$ and $M_N$ have one leaf $L = \{31/54, 49/54 \}$ in common for every $N$. Call $L_N = \{ 1/9 - \delta_N, 4/9 - \delta_N \}$ and $L_\infty = \{1/9, 4/9 \}$ the other leaf in $M_N$ and $M$ respectively. Note the endpoints of $L_N$ and $L_\infty$ differ by $\delta_N$.

Fix $x$ and $y$ in $\mathbb{R}/\mathbb{Z}$ and consider a shortest path between them in $\gamma(M)$. Note this path uses the vertex from identifying the points on the leaf $L_\infty$ at most once. If the path never passes from one endpoint of $L_\infty$ to the other using the vertex induced by identifying the points on $L_\infty$ then the same path exists in $\gamma(M_N)$. If the path does use this identification then we can follow the path from $x$ till we reach the endpoint of $L_\infty$ in $\gamma(M_N)$, then follow a path of length $\delta_N$ to get to the vertex of $\gamma(M_N)$ induced by collapsing $L_N$ to a point, then follow a path of length $\delta_N$ back to the other endpoint of $L_\infty$, and then continue following our path from there. This gives a path in $\gamma(M_N)$ whose length is at most $2\delta_N$ more than the path in $\gamma(M)$.

Similarly we can replace a shortest path in $\gamma(M_N)$ with a path of length at most $2\delta_N$ more in $\gamma(M)$. Thus by definition we have:
\[
md(M_N, M) \leqslant 2\delta_N.
\]
\end{proof}

Here we give a lower bound on the entropy of each $M_N$. Since we do not need to compute the exact entropy we can selectively choose only certain vertices and study their growth rates to get this lower bound and make use of the common form of the leaves of $M_N$. Still the calculations are quite technical.

\begin{lemma}
\label{cubicgrowthlemma}
For each $N \geqslant 2$ let $\lambda_N$ be the unique real number $> 2$ such that:
\[
\lambda_N^N(\lambda_N - 2)^2 = 1.
\]
Then $h(M_N) \geqslant \log \lambda_N$.
\end{lemma}
\begin{proof}
Let $a_N = \frac{1}{3} - \frac{1}{2 \cdot 3^N}$, $b = 13/18$. Note that the endpoints of the leaf $L_N = \{1/9 - \delta_N, 4/9 - \delta_N \}$ map to $a_N$ and the endpoints of the leaf $L = \{ 31/54, 49/54 \}$ map to $b$ under $T(x) = 3x \mod 1$.

Notice that $T^j(a_N) = 1 - \frac{1}{2\cdot 3^{N-j}} \mod 1$ for $1 \leqslant j \leqslant N$ and $T^{N+j}(a_N) = 1/2$ for $j \geqslant 1$. Also we have that $T^j(a_N) \geqslant 17/18 > 49/54$ for all $1 \leqslant j \leqslant N-2$.

\begin{figure}[htbp]
\centering
\begin{tikzpicture}[font=\small]
  \def\R{3.2}
  \pgfmathsetmacro{\leafradius}{\R*sqrt(3)}
  \pgfmathsetmacro{\firstangle}{3220/81}
  \pgfmathsetmacro{\secondangle}{620/3}

  \draw[gray] (0,0) circle[radius=\R cm];

  \draw[thick] (\firstangle:\R)
    arc[start angle={\firstangle-90},
        end angle={\firstangle-150},radius=\leafradius cm]
    node[midway,below=3pt] {$L_5$};

  \draw[thick] (\secondangle:\R)
    arc[start angle={\secondangle-90},
        end angle={\secondangle-150},radius=\leafradius cm]
    node[midway,above=3pt] {$L$};

  \foreach \angle/\value in {
    {3220/81}/{\frac{161}{1458}},
    {12940/81}/{\frac{647}{1458}},
    {620/3}/{\frac{31}{54}},
    {980/3}/{\frac{49}{54}}
  } {
    \node at (\angle:3.7) {$\value$};
  }

  \draw[gray,thin]
    ({3220/9}:\R) -- (3.65,0.15) -- (4,0.15)
    node[right,text=black] {$\frac{161}{162}$};

  \draw[gray,thin]
    ({1060/3}:\R) -- (3.65,-0.55) -- (4,-0.55)
    node[right,text=black] {$\frac{53}{54}$};

  \foreach \angle/\value in {
    60/{\frac16},
    {3220/27}/{a_5=\frac{161}{486}},
    180/{\frac12},
    260/{\frac{13}{18}},
    300/{\frac56},
    340/{\frac{17}{18}}
  } {
    \node at (\angle:3.75) {$\value$};
  }

  \foreach \angle in {
    60,{3220/27},180,260,300,340,{1060/3},{3220/9}
  } {
    \fill (\angle:\R) circle[radius=1.8pt];
  }
\end{tikzpicture}
\caption{$M_5$ and its postcritical set.}
\label{fig:cubic-M5}
\end{figure}
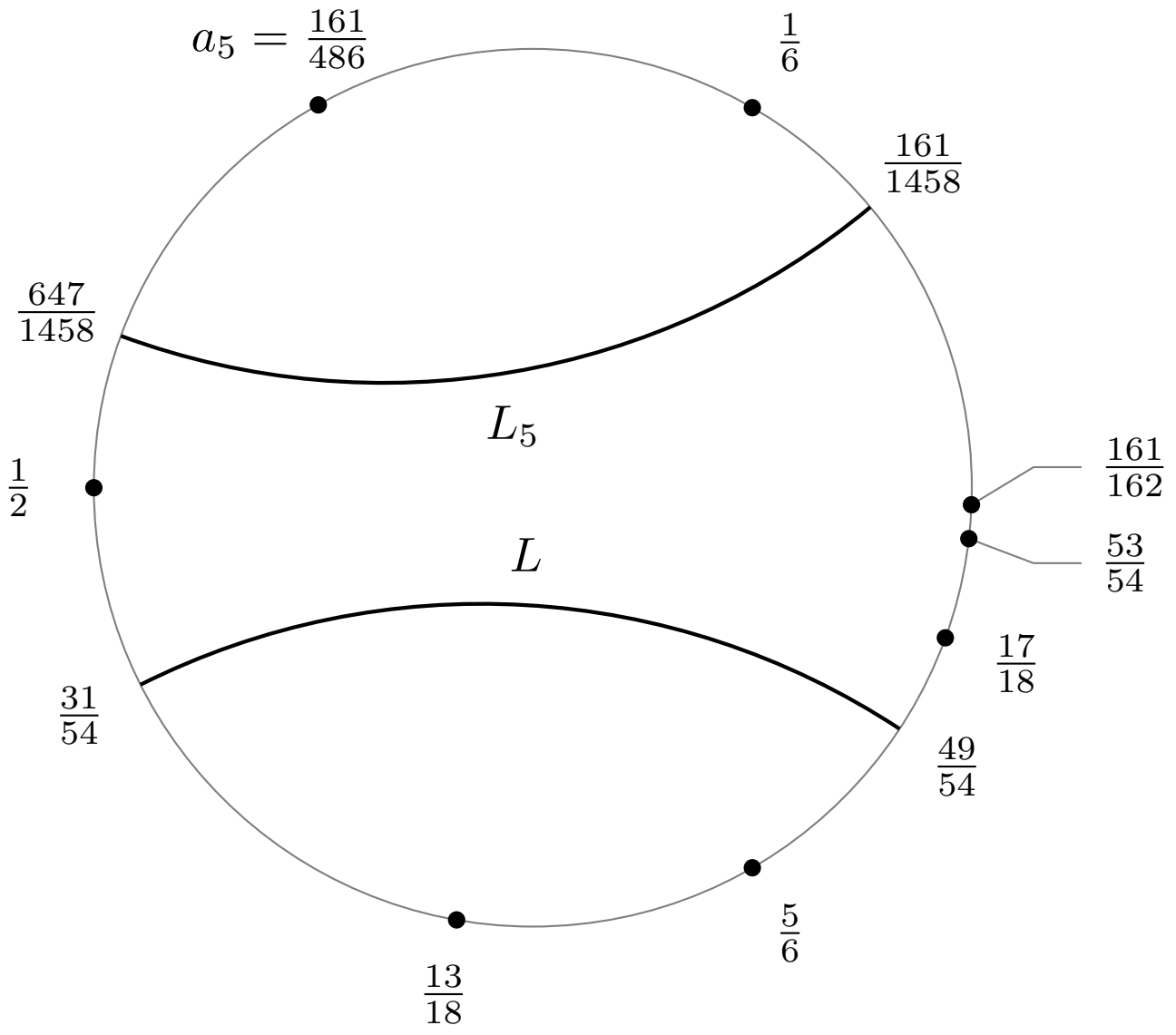

Notice that $L_N$ separates $T^j(a_N)$ from $a_N$ for each $1 \leqslant j \leqslant N-2$ and that $L$ does not. Thus if we let $R_j = (a_N, T^j(a_N))$ then for each $1 \leqslant j \leqslant N-2$ we have edges:
\[
R_j \rightarrow R_1, R_{j+1}.
\]
Since $T^{N-1}(a_N) = 5/6$ we have that $a_N$ and $T^{N-1}(a_N)$ are separated by both leaves. Thus if we let $V = (a_N, 13/18)$ then we have edges:
\[
R_{N-1} \rightarrow R_1, V.
\]
Consider also the vertex $W = (1/6, 13/18)$. Since $a_N$ is separated from $13/18$ by both leaves we have edges:
\[
V \rightarrow R_1, V, W
\]
and:
\[
W \rightarrow V, W.
\]
See figure \ref{fig:cubic-M5} for the $N=5$ picture with the critical orbits drawn. Note there are also edges from $R_{N-1} \rightarrow (13/18, 1/2)$ and $W \rightarrow (a_N, 1/2)$ however our argument will not use these. 

Consider the subgraph $G$ of $G_{M_N}$ consisting of the vertices $R_1, \dots, R_{N-1}$, $V$, $W$, along with the edges between them described above. See figure \ref{fig:cubic-subgraph-M5} for the subgraph when $N=5$.

\begin{figure}[htbp]
\centering
\begin{tikzpicture}[
  >={Stealth[length=1.8mm,width=1.25mm]},
  vertex/.style={
    circle,draw,minimum size=8mm,inner sep=0pt
  },
  every edge/.style={draw,->,semithick},
  font=\normalsize
]
  \node[vertex] (R1) at (0,0) {$R_1$};
  \node[vertex] (R2) at (2,0) {$R_2$};
  \node[vertex] (R3) at (4,0) {$R_3$};
  \node[vertex] (R4) at (6,0) {$R_4$};
  \node[vertex] (V)  at (8,0) {$V$};
  \node[vertex] (W)  at (8,-2) {$W$};

  \path
    (R1) edge (R2)
    (R2) edge (R3)
    (R3) edge (R4)
    (R4) edge (V);

  \draw[->,semithick] (R2.north)
    .. controls (1.5,1.1) and (0.9,1.1)
    .. (R1.30);

  \draw[->,semithick] (R3.north)
    .. controls (3,1.85) and (1,1.85)
    .. (R1.65);

  \draw[->,semithick] (R4.north)
    .. controls (4.4,2.55) and (0.65,2.55)
    .. (R1.100);

  \draw[->,semithick] (V.north)
    .. controls (6,3.25) and (0,3.25)
    .. (R1.135);

  \path
    (V) edge[bend left=22] (W)
    (W) edge[bend left=22] (V);

  \path
    (R1) edge[loop below,min distance=9mm] (R1)
    (V)  edge[loop right,min distance=9mm] (V)
    (W)  edge[loop right,min distance=9mm] (W);
\end{tikzpicture}
\caption{The subgraph $G$ for the $N=5$ case.}
\label{fig:cubic-subgraph-M5}
\end{figure}
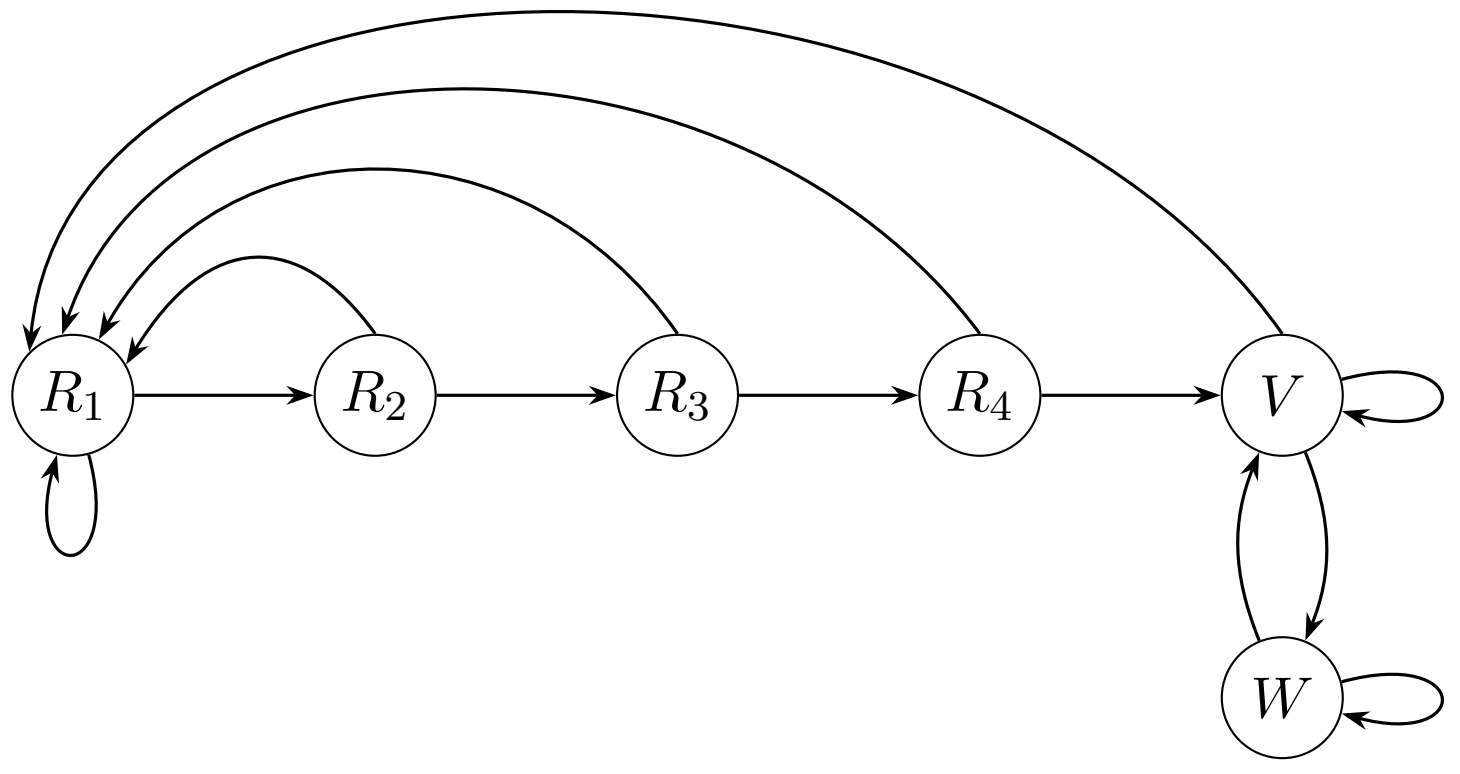

Consider $\lambda$ the leading eigenvalue of $G$. Note $\lambda > 0$ as there is an edge $R_1 \rightarrow R_1$. Moreover, as any vertex in this subgraph has a path to $R_1$ and $R_1$ has a path to any other vertex, the subgraph is strongly connected, hence it has a strictly positive right eigenvector. Let $y$ be such an eigenvector and assume $y_{R_1} = 1$. Then the edges for $V$ and $W$ give:
\[
(\lambda-1)y_V = 1 + y_W, \quad (\lambda - 1)y_W = y_V.
\]
Combining these and simplifying gives:
\[
\lambda(\lambda-2)y_W= 1.
\]
Notice that since $y_W > 0$ we have $\lambda > 2$. For $j \leqslant N-2$ the edges of $R_j$ give:
\[
y_{R_{j+1}} = \lambda y_{R_j} - 1.
\]
Since we chose $y_{R_1} = 1$, this implies that for $j \leqslant N-1$:
\[
y_{R_j} = \frac{(\lambda-2)\lambda^{j-1} + 1}{\lambda - 1}.
\]
At $R_{N-1}$ we have:
\[
\lambda y_{R_{N-1}} = 1 + y_V.
\]
Thus:
\[
(\lambda-2)\lambda^{N-1} + \lambda = (\lambda-1)(1+y_V) = \lambda + y_W.
\]
and so $(\lambda-2)\lambda^{N-1} = y_W$. Since we showed above that $\lambda(\lambda-2)y_W = 1$ this implies:
\[
\lambda^N(\lambda-2)^2 = 1.
\]
Since $G$ is a subgraph of $G_{M_N}$ we have $\lambda(G_{M_N}) \geqslant \lambda$ which completes the proof.
\end{proof}

With these estimates we are finally ready to give a counterexample to the general conjecture.
\begin{proposition}
For $M$ as above we have:
\[
\alpha_{loc}(h, M) \leqslant \frac{\log 2}{2 \log 3} < \frac{h(M)}{\log 3}.
\]
\end{proposition}

\begin{proof}
First note that $h(M) = \log(2)$. Suppose that for some:
\[
\gamma > \frac{\log 2}{2 \log 3},
\]
we had $C, \epsilon > 0$ such that for all $md(M', M) < \epsilon$ we had:
\[
|h(M') - h(M)| \leqslant C md(M', M)^\gamma.
\]
By lemma \ref{cubicmetriclemma} we know we can take $N$ large enough so that $md(M_N, M) \leqslant 3^{-(N+1)} < \epsilon$.

By lemma \ref{cubicgrowthlemma} we know $h(M_N) \geqslant \log \lambda_N$ where $\lambda_N$ is the unique positive number greater than $2$ solving $\lambda^N(\lambda - 2)^2 = 1$.  Then we have:
\[
\log \frac{\lambda_N}{2} \leqslant h(M_N) - h(M) \leqslant C3^{-\gamma(N+1)}.
\]
Since $\log x \geqslant 1-1/x$ for $x > 0$, the definition of $\lambda_N$ lets us estimate the lower bound:
\[
\log \frac{\lambda_N}{2} \geqslant \frac{\lambda_N-2}{\lambda_N} = \lambda_N^{-\frac{(N+2)}{2}}.
\]
Combining inequalities and taking an $N$-th root gives:
\[
\lambda_N^{-\frac{(N+2)}{2N}} \leqslant C^{1/N} 3^{-\gamma\frac{N+1}{N}}.
\]
Since $\lambda_N \rightarrow 2$, taking the limit gives:
\[
2^{\frac{-1}{2}} \leqslant 3^{-\gamma}.
\]
This contradicts $\gamma > \frac{\log2}{2 \log 3}$. Thus:
\[
\alpha_{loc}(h, M) \leqslant \frac{\log 2}{2 \log 3},
\]
as wanted.

\end{proof}

\end{document}